%% file: Paper.tex
\documentclass[preprint,12pt]{article}

\usepackage{amssymb}
\usepackage{amsmath}
\usepackage{moreverb,url}
\usepackage{natbib}
\setcitestyle{square,numbers}
\usepackage{cleveref}
\usepackage{soulutf8}
\usepackage{authblk} 
\usepackage{bigints}
\usepackage{graphicx}
\usepackage{pgfplots}
\usepackage{algorithm, algorithmic}

\renewcommand{\vec}[1]{\textbf{#1}}
\newcommand{\vecx}{\vec{x}} 
\newcommand{\vecv}{\vec{v}} 
\newcommand{\vecxv}{\vecx,\vecv}
\newcommand{\vecxvGrid}[2]{\vecx_{\vec{#1}},\vecv_{\vec{#2}}}
\newcommand{\nablav}{\nabla_\vec{v}} 
\newcommand{\nablax}{\nabla_\vec{x}} 
\newcommand{\supidx}[1]{[#1]}

\newcommand{\parac}[1]{\ensuremath{#1_{z}}}
\newcommand{\perpc}[1]{\ensuremath{#1_{\perp}}}

\title{Convergence of splitting methods on rotating grids for the magnetized Vlasov equation}

\author[1]{Nils Schild}
\author[1]{Mario Räth}
\author[1]{Klaus Hallatschek}
\author[2]{Katharina Kormann}
\affil[1]{Max Planck Institut for Plasma Physics, Germany}
\affil[2]{Ruhr University Bochum, Germany}
\date{} 
\begin{document}
\maketitle

\begin{abstract}
Semi-Lagrangian solvers for the Vlasov system offer noiseless solutions compared
to Lagrangian particle methods and can handle larger time steps compared to
Eulerian methods. In order to reduce the computational complexity of the
interpolation steps, it is common to use a directional splitting. However, this
typically yields the wrong angular velocity. In this paper, we analyze a
semi-Lagrangian method that treats the $\vecv \times \vec B$ term with a
rotational grid and combines this with a directional splitting for the remaining
terms. We analyze the convergence properties of the scheme both analytically and
numerically. The favorable numerical properties of the rotating grid solution are
demonstrated for the case of ion Bernstein waves.
\end{abstract}

%

\section{Introduction and Problem Description}
A kinetic description  of (collisionless) plasmas evolves the phase-space
distribution function $f_s$ of particles of species $s$ (of charge $q_s$ and mass $m_s$ in external and
self-consistent electromagnetic fields $\vec{E}$ and $\vec{B}$ by the so-called
Vlasov equation
\begin{align}
    \partial_t f_s(\vecxv,t) + \vecv\cdot \nablax f_s(\vecxv,t) + \dfrac{q_s}{m_s}(\vec{E}(\vecx,t) + \vecv \times \vec{B})\cdot\nablav f_s(\vecxv,t) = 0.
\end{align}
While this hyperbolic conservation law appears linear for given electromagnetic
fields, the equation is non-linearly coupled to Maxwell's equations for
self-consistent fields. In this paper, we consider a simplified model where the
self-consistent magnetic field is neglected, and the background magnetic field
$\vec{B}_0$ is considered to be a constant field in both time and space-aligned
with the $\hat{e}_z$ axis. Our physical use cases are motivated by ion transport
properties in plasmas. Those can be described by assuming adiabatic electrons,
which handles the electron movement implicitly and a quasi-neutrality
assumption. We normalize physical quantities ($q=m=T=1$) such that the ion
motion in our model is described by
\begin{align}
    \label{num:equ:VlasovEquation}
    \partial_t f(\vecxv,t) &+ \vecv \cdot \nablax f(\vecxv,t) + (\vec E(\vecx,t) + \vecv \times \vec B_0)\cdot \nablav f(\vecxv,t) = 0\\
    \label{num:equ:QuasiNeutralWithAdiabaticElectrons}
    \phi(\vecx,t) &= n(\vecx,t) = \int f(\vecxv,t) \mathrm{d}\vecv, \quad \vec E(\vecx,t) = -\nablax\phi(\vecx,t).
\end{align}
and was used to investigate the limits of gyrokinetics \cite{RaethPhysPlasmas24}
and study turbulence phenomena that are not covered by gyrokinetic models
\cite{RaethPhysRevLett24}. The authors emphasize that the introduced rotating
grid is not limited to this model but can be utilized for any electrostatic 
model with a constant magnetic background field. An example would be the Vlasov
Poisson model is widely used to verify numerical methods. 

The backward semi-Lagrangian method discretizes the distribution function on a
grid. The point-wise solution at a given grid point is propagated forward in
time in two steps: First, the characteristic equations of motion associated with
the hyperbolic conservation properties of \cref{num:equ:VlasovEquation} are
solved backward in time until the previous time step. Then, the solution at the
grid point at the new time step is given by the solution at the previous time
step at this foot of the characteristic curve.  In order to simplify the
solution of the characteristic, it is common to use a directional splitting that
solves the characteristics along one dimension at a time and then combines the
six steps in Lie, Strang, or higher-order splitting method
(cf.~\cite{Cheng1976}). Due to the splitting, the characteristic equations are
not solved exactly and, in particular an inexact solution of the rotation
induced by the term $\vec{v} \times \vec{B}$ can yield a numerical heating of
the plasma as explained in \cite{Schmitz2006}. Therefore, several approaches
have been proposed in the literature to better approximate in particular that
rotation induced by the $\vec{v} \times \vec{B}$ term.

This paper starts with the idea of a rotating grid, as proposed by Kormann,
Reuter \& Rampp \cite{Kormann2019} which removes the rotation induced by
$\vec{v} \times \vec{B}$ from the advection step by pushing it into a coordinate
transform of the computational grid. We extend the previous work by the
transformation of the full Vlasov equation into the computational domain.
Additionally, we provide a convergence analysis for splitting methods on
the rotating grid. 

The rotating grid has two main advantages compared to the direct solution of
\cref{num:equ:VlasovEquation}. Firstly, the rotation removes the $\vecv$
dependence of the $\vecv$-advective part of the Vlasov equation. The remaining
$\vecv$ advection can be split less expensively without loss of accuracy as will
be explained in \cref{sec:SplittingMethods} on splitting methods.  Secondly, the
rotating grid gives more accurate results compared to the splitting schemes
applied to \cref{num:equ:VlasovEquation} with the $\vecv \times \vec B$ term as
is shown in \cref{sec:NumericalResultsWAndWoRotatingGrid}. In certain situations
with strong background fields, the rotational grid can also lead to more local
interpolation stencils, which is of interest in a distributed solution.

Alternative approaches to accurately solve the rotational motion have been
proposed in the literature. Schmitz \& Grauer \cite{Schmitz2006} proposed a
backsubstitution method applied to the Boris scheme.  Bernier, Casas \&
Crouseilles \cite{Bernier2020} propose to decompose a two-dimensional rotational
motion into a product of three shear transformations that amount to
one-dimensional advection steps each. Compared to the latter approach, the use
of a rotational grid has the advantage that the number of split steps is
smaller, which significantly reduces the computational cost to calculate the
solution.

The main goal of this paper is to provide a complete analysis of the
convergence properties of the semi-Lagrangian method with a rotational grid and
a directional splitting.  Convergence of semi-Lagrangian schemes has been
studied in \cite{Besse2008} for the one-dimensional Vlasov--Poisson system and
by Einkemmer \& Ostermann \cite{Einkemmer2014} with a particular focus on the
directional splitting time accuracy. Our analysis extends on the analysis
provided by Einkemmer \& Ostermann for the semi-Lagrangian method without a
rotational grid and builds on the techniques summarized in \cite{Hairer2006} for
the analysis of splitting methods and a Lagrangian-Eulerian viewpoint on the
Vlasov equation \cite{donea2004,Huang2011,Mukhamet2023}. Moreover, we will show for the
example of ion Bernstein waves that solution of superior quality---in
particular with respect to higher modes in both space and time---can be achieved
with the rotational grid compared to a pure directional splitting.

The remainder of the article is organized as follows: In the following section,
we derive the Vlasov equation in the rotational domain and briefly recapture the
semi-Lagrangian method which is applied in the rotating frame.  Section
\ref{sec:SplittingMethods} considers the temporal splitting method and an
analysis of its convergence properties. Numerical results that verify the error
analysis are presented in Section \ref{sec:NumericalResultsWAndWoRotatingGrid}
along with the physical test case of nonlinear ion Bernstein waves demonstrating
the positive effect of the use of the rotational grid semi-Lagrangian method.

\section{Coordinate transformation and semi-Lagrangian discretization}
In this section we first transform the Vlasov equation into the rotating frame.
Afterwards, we briefly introduce the semi-Lagrangian method which is used to 
solve the Vlasov equation for a given initial condition.

\subsection{Rotating velocity frame}
The coordinate transformation of this subsection will remove the $\vecv \times
\vec{B}_0$ term in \cref{num:equ:VlasovEquation} using a rotating velocity grid.

The required moving mesh is based on the coordinate transformation described by
Huang \& Russell \cite[Chap.~3.1]{Huang2011}. In order to derive the rotating
grid we only need to consider the rotational part of
\cref{num:equ:VlasovEquation}
\begin{align}
    \label{num:equ:VlasovEquationRotationPart}
    \partial_t f(\vecv,t) + ( \vecv \times \boldsymbol{\omega}_c)\cdot\nablav f(\vecv,t) = 0
\end{align}
with the cyclotron frequency $\boldsymbol{\omega}_c=q/m\vec{B}_0$. Here and in the following we omit the index $s$ for notational simplicity. The
distribution function shall now be mapped onto a computational domain
$\Omega_\text{C}$ which rotates with respect to the physical velocity domain
$\Omega$. A mapping with the following structure has to be constructed
\begin{align}
    \label{num:equ:MappingMovingMeshToCompDomain}
    \vecv = \vecv(\tilde{\vecv},\tau) : \Omega_\text{C}\times [0,T]  \rightarrow \Omega
\end{align}
where $\tilde{\vecv}$ is the velocity coordinate on the rotating grid and the
time of the rotating grid is the same as on the physical grid $t=\tau$. This
mapping shall remove the rotational part in
\cref{num:equ:VlasovEquationRotationPart}.  First, the derivatives with respect
to $\vecv$ and $t$ are substituted by the derivatives with respect to
$\tilde{\vecv}$ and $\tau$.  The gradient operator with respect to
$\tilde{\vecv}$ can be obtained through the chain rule
\begin{align}
    \label{num:equ:GradientMovingGrid}
    \nablav = \sum_{i}\left(\nablav \tilde{v}_i(\vecv,t)\vert_t\right) \partial_{\tilde{v}_i} = (\vec{J}^{-1})^T \nabla_{\tilde{\vecv}}
\end{align}
where $\tilde{v}_i$ is the $i$-th component of $\tilde{\vecv}$ and
$\vec{J}^{-1}=\dfrac{\partial \tilde{\vecv}}{\partial \vecv}$ is the Jacobian of
the inverted mapping of \cref{num:equ:MappingMovingMeshToCompDomain}.
Additionally, the partial time derivative with respect to $t$ has to be replaced
by a partial derivative with respect to $\tau$ and $\tilde{\vecv}$.
\begin{align}
\label{num:equ:TimeDerivativeMovingGrid}
    \partial_t f\vert_\vecv &= (\partial_\tau f\vert_{\tilde{\vecv}}) (\partial_t \tau) + \sum_{i} (\partial_{\tilde{v}_i} f\vert_\tau) (\partial_t \tilde{v}_i)\\
                            &= (\partial_\tau f\vert_{\tilde{\vecv}}) + (\nabla_{\tilde{\vecv}}f\vert_\tau) \cdot (\partial_t \tilde{\vecv})
\end{align}
Now we insert \cref{num:equ:GradientMovingGrid} and
\cref{num:equ:TimeDerivativeMovingGrid} into the rotational part of the Vlasov
equation \cref{num:equ:VlasovEquationRotationPart}
\begin{align}
    0=&\partial_\tau f(\tilde{\vecv},\tau) +  (\partial_t \tilde{\vecv}) \cdot \nabla_{\tilde{\vecv}}f(\tilde{\vecv},\tau) + (\vecv(\tilde{\vecv},\tau)\times\boldsymbol{\omega}_c) \cdot \left((\vec{J}^{-1})^T \nabla_{\tilde{\vecv}} f(\tilde{\vecv},\tau))\right)\\
    =&\partial_\tau f(\tilde{\vecv},\tau) +  (\partial_t \tilde{\vecv}) \cdot \nabla_{\tilde{\vecv}}f(\tilde{\vecv},\tau) + \left(\vec{J}^{-1}(\vecv(\tilde{\vecv},\tau)\times\boldsymbol{\omega}_c)\right) \cdot \left(\nabla_{\tilde{\vecv}} f(\tilde{\vecv},\tau))\right)
\end{align}
The mapping in \cref{num:equ:MappingMovingMeshToCompDomain} will remove the
rotational part in the Vlasov equation if the following condition is met
\begin{align}
    \partial_t \tilde{\vecv} = -\left(\vec{J}^{-1}(\vecv(\tilde{\vecv},\tau)\times\boldsymbol{\omega}_c)\right).
\end{align}
If we consider a constant background magnetic field in z-direction
$\vec{B}_0=B_0\hat{e}_z$ We can write the equation as
\begin{align}
    \begin{pmatrix}
        \partial_t \tilde{v}_x\\
        \partial_t \tilde{v}_y\\
        \partial_t \tilde{v}_z\\
    \end{pmatrix}
    = 
    \begin{pmatrix}
        - (\partial_{v_x}\tilde{v}_x) \omega_c v_y + (\partial_{v_y}\tilde{v}_x) \omega_c v_x \\
        - (\partial_{v_x}\tilde{v}_y) \omega_c v_y + (\partial_{v_y}\tilde{v}_y) \omega_c v_x \\
        - (\partial_{v_x}\tilde{v}_z) \omega_c v_y + (\partial_{v_y}\tilde{v}_z) \omega_c v_x \\
    \end{pmatrix}
\end{align}
which is satisfied by
\begin{align}
    \label{num:equ:PhysicalDomainToRotatingFrame}
    \tilde{\vecv} &= \vec{D}_{\omega_c}(t)\vecv\\
    \vec{D}_{\omega_c}(t) &= 
    \begin{pmatrix}
        \cos(\omega_c t) & -\sin(\omega_c t) & 0 \\
        \sin(\omega_c t) & \cos(\omega_c t) & 0 \\
        0 & 0 & 1 
    \end{pmatrix}.
\end{align}
It can be verified by insertion and defines the mapping between the
computational and the physical domain in
\cref{num:equ:MappingMovingMeshToCompDomain}. Solving
\cref{num:equ:VlasovEquationRotationPart} in on the rotating domain reduces to a
trivial problem
\begin{align}
    \label{num:equ:VlasovEquationRotationPartCompDomain}
    \partial_\tau f(\tilde{\vecv},\tau) = 0.
\end{align}
Finally, we map the Vlasov equation in \cref{num:equ:VlasovEquation} into the
computational domain substituting $\vecv$ by $\tilde{\vecv}$ and again using
\cref{num:equ:GradientMovingGrid} which results in
\begin{align}
    \label{num:equ:VlasovEquationRotatingFrame}
    \partial_\tau f(\vecx,\tilde{\vecv},\tau) + (\vec{D}_{\omega_c}^{-1}(\tau)\tilde{\vecv})\cdot \nablax f(\vecx,\tilde{\vecv},\tau) + \left(\vec{D}_{\omega_c}(\tau)\dfrac{q}{m}\vec{E}(\vecx,\tau)\right)\cdot \nabla_{\tilde{\vecv}} f(\vecx,\tilde{\vecv},\tau) = 0.
\end{align}
The notation will be simplified in the remainder of the paper. The tilde and $\tau$
will no longer be used to highlight the moving velocity mesh. Only if it is
essential to distinguish between the physical and the computational domain we
will explicitly use $(\tilde{\vecv},\tau)$ instead of $(\vecv,t)$.

\subsection{Solving the Vlasov equation using semi-Lagrangian methods}
\label{subsec:SemiLagrangianMethod}
Before we consider the actual integration methods in the next section we recapture
the basic idea of the semi-Lagrangian method which is our chosen numerical method to
implement the integrators. A detailed discussion on the semi-Lagrangian method
can be found in \cite{Falcone2014}.

The semi-Lagrangian method propagates the distribution function based on the
conservation properties of the hyperbolic partial differential equation. The distribution function $f$ is
conserved along the trajectories of the so-called characteristic curves. The
characteristic curves of the Vlasov equation in \cref{num:equ:VlasovEquation}
are defined by
\begin{align}
    \label{num:equ:VlasovEquationCharacteristics}
    \dfrac{\mathrm d}{\mathrm d t} 
    \begin{pmatrix}
    \vec X(t)\\
    \vec V(t)   
    \end{pmatrix} =
    \begin{pmatrix}
        \vec V(t)\\
        \dfrac{q}{m}( \vec E(\vec X(t)) + \vec V(t) \times \vec B_0)
    \end{pmatrix}.
\end{align}
In the semi-Lagrangian method these characteristics have to be integrated in
time using a phase space grid point $(\vecxvGrid{i}{j})$ as an initial condition 
where $\vec i, \vec j \in \mathrm{N}^d$ are multi-indexes indicating grid
points. We can then use the hyperbolic conservation law to trace the
distribution function after a time step $h$ back to an initial condition
$f_0(\vecxv)$
\begin{align}
    \label{num:equ:BacktracingValue}
    f(\vecxvGrid{i}{j},h) = f_0(\vec X(0;\vecxvGrid{i}{j},h), \vec V(0;\vecxvGrid{i}{j},h)), 
\end{align}
where we denote by $(\vec X(0;\vecx_i,\vecv_j,h),\vec V(0;\vecx_i,\vecv_j,h)$ the solution at time $0$ of the characteristic equations starting at $(\vecx_i,\vecv_j)$ at time $h$ and solved backwards in time. 
The point $(\vec X(0;\vecxvGrid{i}{j},h), \vec V(0;\vecxvGrid{i}{j},h))$ is
usually not a grid point of the initial condition. Therefore, the point has to
be approximated by numerical interpolation
\begin{align}
    f(\vecxvGrid{i}{j},h&) = \\ 
    &I[\{f_0(\vecxvGrid{i}{j})\}]\left( 
        \vecx_{\vec i} + \int_{h}^{0} \vec V(t) \mathrm{d}t,
        \vecv_{\vec j} + \int_{h}^{0} \dfrac{q}{m}(\vec E(\vec X(t)) + \vec V(t) \times \vec B_0)\mathrm{d}t\right).
\end{align}
Here we denoted by $I[\{f_0(\vecxvGrid{i}{j})\}]$ an arbitrary interpolation
procedure that defines an interpolant based on the tuples
$\{((\vecxvGrid{i}{j}),f_0(\vecxvGrid{i}{j}))\}$.

The characteristics of the Vlasov equation in the rotating frame
\cref{num:equ:VlasovEquationRotatingFrame} are given by
\begin{align}
    \label{num:equ:VlasovEquationCharacteristicsRotatingFrame}
    \dfrac{\mathrm d}{\mathrm d t} 
    \begin{pmatrix}
    \vec X(t)\\
    \vec V(t)   
    \end{pmatrix} =
    \begin{pmatrix}
        \vec D^{-1}_{\omega_c}(t)\vec V(t)\\
        \dfrac{q}{m}\vec D_{\omega_c}(t)\vec E(\vec X(t))
    \end{pmatrix}.
\end{align}
such that the distribution function is advected using
\begin{align}
    f(\vecxvGrid{i}{j},h) = 
    I[\{f_0(\vecxvGrid{i}{j})\}]\left( 
        \vecx_{\vec i} + \int_{h}^{0} \vec D^{-1}_{\omega_c}(t)\vec V(t) \mathrm{d}t,
        \vecv_{\vec j} + \int_{h}^{0} \dfrac{q}{m}\vec D_{\omega_c}(t)\vec E(\vec X(t))\mathrm{d}t\right)
\end{align}
Since the electric field is dependent of $f(\vecxv,t)$ thought the field equations, the advection equation is
nonlinear. The numerical analysis of nonlinear equations becomes significantly
more difficult.  Even if the electric field would be simply a constant background field, we
can integrate the characteristic equations but still have to execute an
interpolation step in up to six dimensions which is computationally expensive.

It is therefore desirable to reduce the dimensionality of a single advection
step to reduce the computational effort and ease the numerical analysis of the
advection method. In the next section we use splitting method to split the six
dimensional problem in multiple lower dimensional problems which can be solved
after one another and are simpler to analyze numerically. The interpolations
which we use within this work are briefly described in
\cref{app:InterpolationMethods}.

\section{Splitting methods applied to the Vlasov equation}
\label{sec:SplittingMethods}
In the previous section we introduced the semi-Lagrangian method and defined the 
characteristics for the two representations of the Vlasov equation with and 
without a rotating velocity grid. In this section we utilize splitting methods
to decompose the single 6-D advection equation into multiple 1-D advection steps
to simplify the solution of the characteristics and the numerical analysis. These are recapitulated in the
following \cref{subsec:SplittingMethodsforOdesAndDifferentialOperators}. In our
considerations on splitting methods we assume that the interpolation error is
small, and the splitting error is the dominant error of the splitting.

\subsection{Splitting methods for differential operators}
\label{subsec:SplittingMethodsforOdesAndDifferentialOperators}
We briefly recapture the ideas of splitting methods in a generalized abstract
concept. For in depth explanations, see Hairer, Lubich \&
Wanner~\cite[Sec.~II.5,Sec.~III.5]{Hairer2006}. First, consider differential
operators of the form 
\begin{align}
    D_i = \sum_j A^{\supidx{i}}_j(\vec y) \partial_{y_j} = \vec A^{\supidx{i}}(\vec y)\cdot \nabla_{\vec y},
\end{align}
with the phase space variable $\vec y \in \mathbb{R}^n$. We can build a
hyperbolic PDE based on these abstract differential operators
\begin{align}
    \label{num:equ:AbstractAdvectionEq}
    (D_1 + D_2) g = 0
\end{align}
where we considered two components. The corresponding characteristics are given
by
\begin{align}
    \dot{\vec y} = \vec A^{\supidx{1}}(\vec y) + \vec A^{\supidx{2}}(\vec y)
\end{align}
with $\dot{\vec y}$ being the time derivative of $\vec y$. We assume that we can
integrate the splitted ordinary differential equations (ODE)
\begin{align}
    \dot{\vec y} = \vec A^{\supidx{1}}(\vec y) \qquad  \dot{\vec y} = \vec A^{\supidx{2}}(\vec y)
\end{align}
exactly and that the solution of the ODE is described by the flow $\vec y =
\varphi_t^{\supidx{i}}(\vec y_0)$ with $i=1,2$.

Then the evolution of any differentiable function $g : \mathbb{R}^n \rightarrow
\mathbb{R}$ in \cref{num:equ:AbstractAdvectionEq} can be approximated using
exponential integrators. Let us assume the flow of the characteristics, can be
separated into two components, which can be integrated exactly. Then we can
advance an initial condition of $g$ in time using $D_i$ with exponential
integrators
\begin{align}
    \label{num:equ:ExponentialIntegratorDerivation}
    g(\varphi^{\supidx{i}}_h(\vec y_0)) &= \left(\sum_{n\geq 0} \dfrac{h^n}{n!} (D_i^n g)\right) (\varphi^{\supidx{i}}_h(\vec y_0)\vert_{h=0})\\
                                                   &= \sum_{n\geq 0} \dfrac{h^n}{n!} (D_i^n g) (\vec y_0)= \exp(h D_i)g(\vec y_0).
\end{align}
Here the previously introduced differential operator $D_i g(\vec y) = \vec
A^{\supidx{i}}(\vec y)\cdot \nabla g(\vec y)$ has been used to substitute the
derivative $\mathrm{d}^n/\mathrm{d}t^n g(\varphi^{\supidx{i}}_t(\vec y_0)) =
(D_i^n g)(\varphi_t^{\supidx{i}}((\vec y_0)))$. 

If we substitute $D_i$ by $D=D_1+D_2$ in
\cref{num:equ:ExponentialIntegratorDerivation}, the exponential integrator for
\cref{num:equ:AbstractAdvectionEq} is given by
\begin{align}
    g(\varphi_h(\vec y_0)) = \exp(h D) g(\vec y_0) = \exp(h (D_1 + D_2)) g(\vec y_0)
\end{align}
We can also apply \cref{num:equ:ExponentialIntegratorDerivation} twice to split
the integration step into two parts and advect $g$ with both parts separately.
Advancing the initial condition $y_0$ firstly by $D_1$ and secondly by $D_2$ we receive
\begin{align}
    \label{num:equ:LieSplittingExponentialIntegrator}
    g((\varphi_h^{\supidx{2}}\circ \varphi_h^{\supidx{1}})(\vec y_0)) = g(\varphi_h^{\supidx{2}}(\varphi_h^{\supidx{1}}(\vec y_0))) = \exp(h D_1) \exp(h D_2) g(\vec y_0),
\end{align}
if we utilize \cref{num:equ:ExponentialIntegratorDerivation} recursively.  The
crucial part to consider here is that if $D_1$ and $D_2$ do not 
commute, such that $[D_1,D_2]\neq 0$, we can not simply merge the two
exponential integrators into one which is equal to $\exp(h(D_1 + D_2))$. The
relation
\begin{align}
    \label{num:equ:TwoStepExpIntegrator}
    \exp(h D_2) \exp(h D_1) = \exp(Z(h,D_1,D_2)) \neq \exp(h(D_1 + D_2))
\end{align}
has to be taken into account. Here $Z(h,D_1,D_2)$ is an expansion of terms in powers of $h$ defined
through the Baker-Campbell-Hausdorff (BCH) formula. The splitted integrator is
exact up to matching orders of the terms in $Z(h,D_1,D_2)$ and the exponent of
the right-hand side. The simple integrator defined by
\cref{num:equ:LieSplittingExponentialIntegrator} is locally of order $O(h^2)$
and called Lie-Splitting. A Strang-Splitting is given by
\begin{align}
  \exp(h/2D_1)\exp(hD_2)\exp(h/2D_1)  
\end{align}
and is of locally of order $O(h^3)$.

Higher order integrators can be obtained by using $m$ steps instead of two
\begin{align}
    \exp(b_m h D_1)\exp(a_m h D_2)\exp(b_{m-1} h D_1)...\exp(a_2 h D_1)\exp(b_1 h D_2)\exp(a_1 h D_1)g(\vec y_0).
\end{align}
The coefficients $a_1,b_1,\dots,a_m,b_m$ have to be determined using
\cref{num:equ:TwoStepExpIntegrator}, such that the terms of the polynomial
$Z(h,D_1,D_2)$ vanish up to a given order $h^p$ to obtain an integrator of
order $p$. A detailed explanation on order conditions through the BCH formula is
given by Hairer, Lubich \& Wanner~\cite[Sec.III.4,III.5]{Hairer2006}.

If the integration is actually implemented through an algorithm, an appropriate
numerical method has to be chosen to actually carry out the integration. This
method might introduce further errors based on the time step and errors based on
the discretization of space. We use the semi-Lagrangian method
\cref{subsec:SemiLagrangianMethod} to explicitly implement the exponential
integrators.

\subsection{Splitting propagator for the Vlasov equation in the physical frame}
\label{subsec:IntegratorsForTheVlasovEquation}
After recapturing the basic ideas of splitting methods we will reduce the 6-D Vlasov 
equation in \cref{num:equ:VlasovEquation} down to
multiple 1-D advection problems.  Additionally, we construct an integrator
which is second order accurate, that is each single time step is  required to be
$O(h^3)$ accurate which is equivalent to a Strang splitting. First we need to identify
the differential operators which can be splitted in the Vlasov equation. The
$\vecx$-advectivion operator is $\vecv\cdot \nablax$ which defines the transport
properties in the spatial domain. The velocity domain transport is defined through the second 
differential operator $(\vec E(\vecx,t) + \vecv \times \vec B_0)\cdot \nablav$.  We
split according to these two operators which gives
\begin{align}
    \vec A_\vecx \cdot \nablax = \vecv \cdot \nablax \qquad \vec A_\vecv\cdot \nablav = (\vec E(\vecx,t) + \vecv \times \vec B_0)\cdot \nablav.
\end{align}
The ODEs defining the splitted flows are then given by
\begin{align}
    \dfrac{\mathrm d}{\mathrm d t} 
    \begin{pmatrix}
    \vec X(t)\\
    \vec V(t)   
    \end{pmatrix} =
    \begin{pmatrix}
        \vecv\\
        0
    \end{pmatrix}
    \qquad 
    \dfrac{\mathrm d}{\mathrm d t} 
    \begin{pmatrix}
    \vec X(t)\\
    \vec V(t)   
    \end{pmatrix} =
    \begin{pmatrix}
        0\\
        (\vec E(\vec X(t)) + \vec V \times \vec B_0)
    \end{pmatrix}
\end{align}
In our physically interesting examples of \cref{subsec:NonlinearVlasovExamples},
a field equation is coupled to the Vlasov equation to determine the electric
field. These field equations depend on $f$ only through the particle density
$n$. Therefore, during the $\vecv$ advection based on $\vec A_\vecv$, $f$
changes only with respect to $\vecv$. The particle density $n$ and therefore
also the electric field $\vec E(\vec X(t))$ do not change such that we can drop
the explicit time dependence of the latter in the substep, and we obtain an
autonomous ODE defining the flow of $f$. For autonomous ODE we can utilize the
previously introduced framework of exponential integrators to propagate the
distribution function in time. A rigorous proof of the second order accuracy of
this splitting was provided in Einkemmer \& Ostermann\cite{Einkemmer2014}.

So far we have reduced the six dimensional Vlasov equation into two 3-D problems
that can be solved to propagate the Vlasov equation using the vector fields
$\vec A_\vecx$ and $\vec A_\vecv$. Further reduction to multiple 1-D problems is
achieved by again splitting the vector fields in the spatial domain $A_{i}$ with
$i=(x,y,z)$ and the velocity domain $A_{j}$ with $j=(v_x,v_y,v_z)$. A second
order integrator based on Strang-Splitting is then given by
\begin{align}
    f(\vecxv,h) + O(h^3) = &\exp\left(\dfrac{h}{2}A_{v_x}\partial_{v_x}\right)\exp\left(\dfrac{h}{2}A_{v_y}\partial_{v_y}\right)\exp\left(\dfrac{h}{2}A_{v_z}\partial_{v_z}\right)\nonumber\\
                  &\exp\left(\dfrac{h}{2}A_{x}\partial_x\right)\exp\left(\dfrac{h}{2}A_{y}\partial_y\right)\exp\left(\dfrac{h}{2}A_{z}\partial_z\right)\nonumber\\
                  &\exp\left(\dfrac{h}{2}A_{z}\partial_z\right)\exp\left(\dfrac{h}{2}A_{y}\partial_y\right)\exp\left(\dfrac{h}{2}A_{x}\partial_x\right)\nonumber\\
                  &\exp\left(\dfrac{h}{2}A_{v_z}\partial_{v_z}\right)\exp\left(\dfrac{h}{2}A_{v_y}\partial_{v_y}\right)\exp\left(\dfrac{h}{2}A_{v_x}\partial_{v_x}\right)f_0(\vecxv)
\end{align}
Since the operators in the spatial domain commute $[A_{x_i},A_{x_j}] = 0$ we can
reduce the computational complexity of the problem by switching and merging
operators working on the same axis into a single operation such that we can
reduce twelve operations to nine 
\begin{align}
    \label{num:equ:IntegratorVlasovEquation}
    f(\vecxv,h) + O(h^3) = &\exp\left(\dfrac{h}{2}A_{v_x}\partial_{v_x}\right)\exp\left(\dfrac{h}{2}A_{v_y}\partial_{v_y}\right)\exp\left(\dfrac{h}{2}A_{v_z}\partial_{v_z}\right)\nonumber\\
                  &\exp\left(hA_{x}\partial_x\right)\exp\left(hA_{y}\partial_y\right)\exp\left(hA_{z}\partial_z\right)\nonumber\\
                  &\exp\left(\dfrac{h}{2}A_{v_z}\partial_{v_z}\right)\exp\left(\dfrac{h}{2}A_{v_y}\partial_{v_y}\right)\exp\left(\dfrac{h}{2}A_{v_x}\partial_{v_x}\right)f_0(\vecxv)
\end{align}
This operator gives a convergence rate of order two in the time discretization
if other errors sources depending on the time step can be neglected.

The subsection is concluded by combining the splitting with the semi-Lagrangian
method of \cref{subsec:SemiLagrangianMethod}. The splitted 
operators only require to solve a 1-D advection problem with a constant
advection coefficient such that the integral solutions of
\cref{num:equ:VlasovEquationCharacteristics} reduce to
\begin{align}
    \label{num:equ:BacktracingTrajectory}
    Y(t_0) = y + \int_{t_0 + h}^{t_0} c \mathrm{d}t = y - h c
\end{align}
with $Y(0)\in(\vec X(0),\vec V(0))$ and $y\in(\vecxvGrid{i}{j})$.
The shifts are explicitly given by.
\begin{align}
    &c_{x}   = v_x\quad c_{y} = v_y\quad c_{z} = v_z\\
    &c_{v_x} = \dfrac{q}{m}E_x(\vecx) + \dfrac{q}{m}v_{y}B_0\quad c_{v_y} = \dfrac{q}{m}E_y(\vecx) - \dfrac{q}{m}v_{x}B_0\quad c_{v_z} = \dfrac{q}{m}E_z(\vecx)
\end{align}
In \cref{num:alg:SemiLagrangianPhysicalDomain} all steps are combined to provide
the solver for the Vlasov equation. The solver will be compared to the solution
on a rotating grid in \cref{subsec:NonlinearVlasovExamples}. The subscript of
the interpolation indicates the direction of the 1-D interpolation.
\begin{algorithm}
\caption{\label{num:alg:SemiLagrangianPhysicalDomain}Solve \cref{num:equ:VlasovEquation} using operator splitting described in \cref{subsec:IntegratorsForTheVlasovEquation} and the semi-Lagrangian method of \cref{subsec:SemiLagrangianMethod}}
\begin{algorithmic}[1]
\STATE{Initial condition $f_0(\vecxv)$, Time step $h$, Final time $t$}
\WHILE{$t_0 < t$}
\STATE{$f^{\supidx{1}}(\vecxvGrid{i}{j}) = I_{v_x}[\{f_0(\vecxvGrid{i}{j})\}]           \left(\vecx_{\vec i},\vecv_{\vec j} - h/2 (c_{v_x},0,0)\right)$}
\STATE{$f^{\supidx{2}}(\vecxvGrid{i}{j}) = I_{v_y}[\{f^{\supidx{1}}(\vecxvGrid{i}{j})\}]\left(\vecx_{\vec i},\vecv_{\vec j} - h/2 (0,c_{v_y},0)\right)$}
\STATE{$f^{\supidx{3}}(\vecxvGrid{i}{j}) = I_{v_z}[\{f^{\supidx{2}}(\vecxvGrid{i}{j})\}]\left(\vecx_{\vec i},\vecv_{\vec j} - h/2 (0,0,c_{v_z})\right)$}
\STATE{$f^{\supidx{4}}(\vecxvGrid{i}{j}) = I_{x}[\{f^{\supidx{3}}(\vecxvGrid{i}{j})\}]\left(\vecx_{\vec i} - h (c_{x},0,0)\right)$}
\STATE{$f^{\supidx{5}}(\vecxvGrid{i}{j}) = I_{y}[\{f^{\supidx{4}}(\vecxvGrid{i}{j})\}]\left(\vecx_{\vec i} - h (0,c_{y},0)\right)$}
\STATE{$f^{\supidx{6}}(\vecxvGrid{i}{j}) = I_{z}[\{f^{\supidx{5}}(\vecxvGrid{i}{j})\}]\left(\vecx_{\vec i} - h (0,0,c_{z})\right)$}
\STATE{$\vec E(\vecx_{\vec i})$          = solve electric field($f^{\supidx{6}}(\vecxvGrid{i}{j})$)}
\STATE{$f^{\supidx{7}}(\vecxvGrid{i}{j}) = I_{v_z}[\{f^{\supidx{6}}(\vecxvGrid{i}{j})\}]\left(\vecx_{\vec i},\vecv_{\vec j} - h/2 (0,0,c_{v_z})\right)$}
\STATE{$f^{\supidx{8}}(\vecxvGrid{i}{j}) = I_{v_y}[\{f^{\supidx{7}}(\vecxvGrid{i}{j})\}]\left(\vecx_{\vec i},\vecv_{\vec j} - h/2 (0,c_{v_y},0)\right)$}
\STATE{$f(\vecxvGrid{i}{j},t_0)          = I_{v_x}[\{f^{\supidx{8}}(\vecxvGrid{i}{j})\}]\left(\vecx_{\vec i},\vecv_{\vec j} - h/2 (c_{v_x},0,0)\right)$}
\STATE{Set $t_0 = t_0 + h$ and $f_0(\vecxvGrid{i}{j})=f(\vecxvGrid{i}{j},t_0)$}
\ENDWHILE
\RETURN $f(\vecxv,t)$
\end{algorithmic}
\end{algorithm}

\subsection{Splitting propagator for the Vlasov equation in the rotating frame}
\label{subsec:IntegratorsForTheVlasovEquationRotatingFrame}
In the last subsection the Vlasov equation has been splitted into multiple 1-D
problems. This subsection focuses on the Vlasov equation in the rotating frame
defined by \cref{num:equ:VlasovEquationRotatingFrame}. We can again identify the
differential operators in the spatial domain and the velocity domain respectively
\begin{align}
    \vec A_\vecx = \vec D_{\omega_c}^{-1}(t)\vecv \cdot \nablax \qquad \vec A_\vecv = \vec D_{\omega_c}(t) \vec E(\vecx,t) \cdot \nablav.
\end{align}
The ODEs defining the splitted flows are then given by
\begin{align}
    \label{num:equ:VlasovEquationCharacteristicsRotatingFrameSplitted}
    \dfrac{\mathrm d}{\mathrm d t} 
    \begin{pmatrix}
    \vec X(t)\\
    \vec V(t)   
    \end{pmatrix} =
    \begin{pmatrix}
        \vec D_{\omega_c}^{-1}(t)\vec V\\
        0
    \end{pmatrix}
    \qquad 
    \dfrac{\mathrm d}{\mathrm d t} 
    \begin{pmatrix}
    \vec X(t)\\
    \vec V(t)   
    \end{pmatrix} =
    \begin{pmatrix}
        0\\
        \vec D_{\omega_c}(t) \vec E(\vec X(t))
    \end{pmatrix}.
\end{align}
We can again follow the arguments of the last subsection to drop the time
dependence of the electric field $\vec E(\vec X(t))$ due to the constant spatial
properties of the distribution function $f$ during the advection step in
the velocity domain. The explicit time dependence of the
rotation matrices $\vec D_{\omega_c}(t)$ and $\vec D_{\omega_c}^{-1}(t)$ can not
be removed from the ODEs such that we do not obtain autonomous ODEs to which we
could apply the framework of exponential integrators. Fortunately, the rotation
matrices are known explicitly and not complex such that order conditions for the
flows can be derived by solving the ODEs and calculate the flows explicitly
which will be done in the following.

The exact flow that needs to be solved in the rotating frame is given by
\cref{num:equ:VlasovEquationCharacteristicsRotatingFrame}. Advancing the initial
condition $(\vecxv)$ by a time step of length $h$ starting from $t_0$ will be
denoted by the mapping
\begin{align}
    \varphi_{t_0 + h, t_0} : (\vecxv) \mapsto (\vec X(t_0 + h), \vec V(t_0 + h)).
\end{align}
The approximated flows defined by the ODEs in 
\cref{num:equ:VlasovEquationCharacteristicsRotatingFrameSplitted} are
superscripted by the coordinates which are advected by the flow map
\begin{align}
    &\varphi_{t_0 + h, t_0}^{\supidx{\vecx}} : (\vecxv) \mapsto \left(\vecx + \int_{t_0}^{t_0+h} \vec D_{\omega_c}(s)\vecv \mathrm{d} s, \vecv\right) \\
    &\varphi_{t_0 + h, t_0}^{\supidx{\vecv}} : (\vecxv) \mapsto \left(\vecx, \vecv + \int_{t_0}^{t_0+h} \vec D_{\omega_c}(s)\vec E(\vecx) \mathrm{d} s \right).
\end{align}
We derive the convergence order of an explicit splitting that approximates the
flow map $\varphi_{t_0+h,t_0}(\vecxv)$ globally up to second order in time
\begin{align}
    \label{num:equ:RotatingGridFlowMapsSplittings}
    \varphi_{t_0+h,t_0}(\vecxv) = \left(\varphi_{t_0 + h, t_0+h/2}^{\supidx{\vecv}} \circ \varphi_{t_0 + h, t_0}^{\supidx{\vecx}} \circ \varphi_{t_0 + h/2, t_0}^{\supidx{\vecv}}\right)(\vecxv) + R(h).
\end{align}
If the residual $R(h)$ only contains components of $O(h^3)$ the splitting has
the desired convergence rate. The global convergence order of $O(h^2)$ can then
be proved by means of standard arguments of consistency and stability.

We start with the exact expression of the flow and transform it into an
expression for \cref{num:equ:RotatingGridFlowMapsSplittings}. Into the exact
expression we inserted the approximations $\overline{\vecx}$ and
$\overline{\vecv}$ to receive a link between the splitted flow map and the exact
flow map. The approximation $\overline{\vecx}$ is the $\vecx$ component of
$(\varphi_{t_0 + h, t_0}^{\supidx{\vecx}}\circ\varphi_{t_0 +
h/2,t_0}^{\supidx{\vecv}})(\vecx,\vecv)$.  The approximation $\overline{\vecv}$
is the $\vecv$ component of $\varphi_{t_0 +
h/2,t_0}^{\supidx{\vecv}}(\vecx,\vecv)$. Reorganizing the obtained components 
provides the splitted flow maps of the ODEs in
\cref{num:equ:RotatingGridFlowMapsSplittings} as well as the residual $R(h)$.
\begin{align}
    \varphi_{t_0+h,t_0}(\vecxv) &= &&
    \begin{pmatrix}
        \vec X(t_0 + h)\\
        \vec V(t_0 + h)
    \end{pmatrix}\\
    &=
    \begin{pmatrix}
        \vecx\\
        \vecv
    \end{pmatrix}  
    && + \bigintsss_{t_0}^{t_0+h}
    \begin{pmatrix}
        \vec D^{-1}_{\omega_c}(s) \vec V(s)\\
        \vec D_{\omega_c}(s)\vec E(\vec X(s))
    \end{pmatrix} \mathrm{d} s\\
    &=
    \begin{pmatrix}
        \vecx\\
        \vecv
    \end{pmatrix} 
    && + \bigintsss_{t_0}^{t_0+h/2}
    \begin{pmatrix}
        0\\
        \vec D_{\omega_c}(s)\vec E(\vec X(s))
    \end{pmatrix} \mathrm{d} s\\ 
    & &&+ \bigintsss_{t_0}^{t_0+h}
    \begin{pmatrix}
        \vec D^{-1}_{\omega_c}(s) \vec V(s)\\
        0
    \end{pmatrix} \mathrm{d} s \nonumber\\
    & && + \bigintsss_{t_0 + h/2}^{t_0+h}
    \begin{pmatrix}
        0\\
        \vec D_{\omega_c}(s)\vec E(\vec X(s))
    \end{pmatrix} \mathrm{d} s\\
    &=
    \begin{pmatrix}
        \vecx\\
        \vecv
    \end{pmatrix} 
    && + \bigintsss_{t_0}^{t_0+h/2}
    \begin{pmatrix}
        0\\
        \vec D_{\omega_c}(s)\left[\vec E(\vec X(s)) + \left(\vec E(\vecx) - \vec E(\vecx)\right)\right]
    \end{pmatrix} \mathrm{d} s \nonumber\\
    & &&+ \bigintsss_{t_0}^{t_0+h}
    \begin{pmatrix}
        \vec D^{-1}_{\omega_c}(s) \left[\vec V(s) + \left(\overline{\vecv} - \overline{\vecv})\right)\right]\\
        0
    \end{pmatrix} \mathrm{d} s \nonumber\\
    & &&+ \bigintsss_{t_0 + h/2}^{t_0+h}
    \begin{pmatrix}
        0\\
        \vec D_{\omega_c}(s) \left[\vec E(\vec X(s)) + (\vec E(\overline{\vecx}) - \vec E(\overline{\vecx}))\right]
    \end{pmatrix} \mathrm{d} s \\
    &= &&(\varphi_{t_0 + h/2, t_0}^{\supidx{\vecv}} \circ \varphi_{t_0 + h, t_0}^{\supidx{\vecx}} \circ \varphi_{t_0 + h, t_0 + h/2}^{\supidx{\vecv}})(\vecxv) + R(h).
\end{align}
We consider three terms of the residual
\begin{align}
    R(h) &= R_1(h) + R_2(h) + R_3(h),
\end{align}
which are given by
\begin{align}
    \label{num:equ:ResiduumV1Step}
    &R_1(h) = \bigintsss_{t_0}^{t_0+h/2}
    \begin{pmatrix}
        0\\
        \vec D_{\omega_c}(s)\left[\vec E(\vec X(s)) - \vec E(\vecx)\right]
    \end{pmatrix} \mathrm{d} s \\
    R_2(h) &= \bigintsss_{t_0}^{t_0+h}
    \begin{pmatrix}
        \vec D^{-1}_{\omega_c}(s) \left[\vec V(s) - \overline{\vecv}\right]\\
        0
    \end{pmatrix} \mathrm{d} s \\
    \label{num:equ:ResiduumXStep}
    &= \bigintsss_{t_0}^{t_0+h}
    \begin{pmatrix}
        \vec D^{-1}_{\omega_c}(s) \left[\vec V(s) - (\vecv + \int_{t_0}^{t_0+h/2} \vec D_{\omega_c}(s') \vec E(\vecx) \mathrm{d} s')\right]\\
        0
    \end{pmatrix} \mathrm{d} s \\
    \label{num:equ:ResiduumV2Step}
    R_3(h) &= \bigintsss_{t_0 + h/2}^{t_0+h}
    \begin{pmatrix}
        0\\
        \vec D_{\omega_c}(s) \left[\vec E(\vec X(s)) - \vec E\left(\overline{\vecx} \right)\right]
    \end{pmatrix} \mathrm{d} s,
\end{align}
where $\overline{\vecx} =  \vecx + \int_{t_0}^{t_0+h}\vec
D_{\omega_c}^{-1}(s')\left(\vecv + \int_{t_0}^{t_0+h/2}\vec
D_{\omega_c}(s'')\vec E(\vecx)\mathrm{d} s''\right) \mathrm{d} s'$.  If the
integrals
\cref{num:equ:ResiduumV1Step,num:equ:ResiduumXStep,num:equ:ResiduumV2Step} only
contain terms of order $O(h^3)$ our splitted flow map has the required
convergence properties. 

We consider only small time steps $h$ such that we can expand the integral
solution of \cref{num:equ:VlasovEquationCharacteristicsRotatingFrame} and remove
higher order terms 
\begin{align}
\label{num:equ:VlasovEquationCharacteristicsRotatingFrameIntegralSolutionX}
\vec X(s) &= \vecx + \int_{t_0}^{s} \vec D_{\omega_c}^{-1}(s) \vec V(s') \mathrm{d}s' = \vecx + \int_{t_0}^{s} (\vec D_{\omega_c}^{-1}(s) \vec v + O(s))\mathrm{d}s'\\
\vec V(s) &= \vecv + \int_{t_0}^{s} \vec D_{\omega_c}(s) \vec E(\vec X(s')) \mathrm{d}s'\\
\label{num:equ:VlasovEquationCharacteristicsRotatingFrameIntegralSolutionV}
&= \vecv + \int_{t_0}^{s} \left(\vec D_{\omega_c}(s') \vec E(\vec x) + \int_{t_0}^{s'}\vec D_{\omega_c}^{-1}(s'')\vecv \mathrm{d} s'' \nablax \vec E(\vecx) + O(s^2)\right)\mathrm{d}s',
\end{align}
where we inserted \cref{num:equ:VlasovEquationCharacteristicsRotatingFrame} into
the second integral to expand the velocity advection.

\paragraph{Residuals $R_1(h)+R_3(h)$} We first estimate a residual for the
$R_1(h)+R_3(h)$. We can insert the integral solution
\cref{num:equ:VlasovEquationCharacteristicsRotatingFrameIntegralSolutionX} into
$R_3(h)$, expand both expressions for the electric field $\vec E$ with regard to
the time shift given by the integral, and keep the terms up to $O(h)$ which is
sufficient to show that the residuum is $O(h^3)$. The intermediate steps are
omitted in the following
\begin{align}
R_3(h) &=&&\bigintsss_{t_0 + h/2}^{t_0+h}\vec D_{\omega_c}(s)  \left[ \vec E\left(\vecx + \int_{t_0}^{s} (\vec D_{\omega_c}^{-1}(s) \vec v + O(h))\mathrm{d}s'\right) - \right.\\
    &                                                           && \left. \vec E\left( \vecx + \int_{t_0}^{t_0+h}\vec D_{\omega_c}^{-1}(s')(\vecv + \int_{t_0}^{t_0+h/2}\vec D_{\omega_c}(s'')\vec E(\vecx)\mathrm{d} s'') \mathrm{d} s'\right)\right]\mathrm{d} s\\
&= &&\bigintsss_{t_0 + h/2}^{t_0+h}\vec D_{\omega_c}(s)  \left[ \left(\int_{t_0}^{s}\vec D_{\omega_c}^{-1}(s') \mathrm{d}s' - \int_{t_0}^{t_0+h}\vec D_{\omega_c}^{-1}(s')\mathrm{d}s'\right) \nablax \vec E(\vecx) \right]\mathrm{d} s\\
& &&+ O(h^3).
\end{align}
The residual $R_1(h)$ can be expanded using
\cref{num:equ:VlasovEquationCharacteristicsRotatingFrameIntegralSolutionX} as
well such that it reduces to
\begin{align}
R_1(h) = \int_{t_0}^{t_0+h/2}\left(\vec D_{\omega_c}(s) \int_{t_0}^{s}\vec D_{\omega_c}^{-1}(s') \mathrm{s'} \nablax \vec E(\vecx) \mathrm{d}s'  \right)\mathrm{d}s+ O(h^3) .
\end{align}
We can sum both residuals and use that the argument of the rotation matrices is
$O(h)$ due to the integral boundaries such that we can expand the matrix and
integrate only the first non-zero component which is the unity matrix
\begin{align}
&\Vert R_1(h) + R_3(h)\Vert \\
&= \left\Vert \left(\int_{t_0}^{t_0+h} \int_{t_0}^{s} \vec D_{\omega_c}(s-s') \mathrm{d}s' \mathrm{d}s - \int_{t_0+h/2}^{t_0+h} \int_{t_0}^{t_0+h} \vec D_{\omega_c}(s-s')\mathrm{d}s' \mathrm{d}s \right) \nablax \vec E(\vecx)\right\Vert\\
&= \left\Vert \left( \int_{t_0}^{t_0+h} \int_{t_0}^{s} (\mathrm{1}+O(h)) \mathrm{d}s' \mathrm{d}s - \int_{t_0+h/2}^{t_0+h} \int_{t_0}^{t_0+h} (\mathrm{1}+O(h)) \mathrm{d}s' \mathrm{d}s \right) \nablax \vec E(\vecx)\right\Vert\\
&= O(h^3).
\end{align}
\paragraph{Residual $R_2(h)$} The second component of the residual can be
considered on its own. We first insert the integral solution for $\vec V(s)$
given by
\cref{num:equ:VlasovEquationCharacteristicsRotatingFrameIntegralSolutionV}.
Afterwards, we again expand the nonlinearity $\vec E(\vec X(s))$ using the
integral
\cref{num:equ:VlasovEquationCharacteristicsRotatingFrameIntegralSolutionX}.
The remaining integral has the same structure as the final integral of the
previous paragraph.
\begin{align}
&R_2(h) \\
&= \left\Vert \bigintsss_{t_0}^{t_0+h}
    \vec D^{-1}_{\omega_c}(s) \left(\int_{t_0}^{s} \vec D_{\omega_c}(s') \vec E(\vec X(s'))\mathrm{d} s' - \int_{t_0}^{t_0+h/2} \vec D_{\omega_c}(s') \mathrm{d} s' \vec E(\vecx) \right)\mathrm{d} s \right\Vert\\
    &=\left\Vert \bigintsss_{t_0}^{t_0+h}
    \vec D^{-1}_{\omega_c}(s) \left(\int_{t_0}^{s} \vec D_{\omega_c}(s') (\vec E(\vecx) + O(h)) \mathrm{d} s' - \int_{t_0}^{t_0+h/2} \vec D_{\omega_c}(s') \mathrm{d} s' \vec E(\vecx) \right)\mathrm{d} s \right\Vert\\
    &= \left\Vert\bigintsss_{t_0}^{t_0+h}
    \left(\int_{t_0}^{s} \vec D_{\omega_c}^{-1}(s-s') \mathrm{d} s' - \int_{t_0}^{t_0+h/2} \vec D_{\omega_c}^{-1}(s-s') \mathrm{d} s'  \right)\mathrm{d} s \vec E(\vecx)\right\Vert + O(h^3)\\
    &= O(h^3).
\end{align}
Therefore, we have discussed that all components of $\Vert R(h)\Vert$ are of
order $O(h^3)$ locally such that we achieve an overall global convergence order
of $O(h^2)$.

The subsection is concluded by merging steps within the integrator to reduce the
computational effort of a step moving from $t_0$ to $t_0 + h$.  A significant
difference to the splitting of the last subsection is that in the rotating frame
also the flows, which define the transport properties of the velocity domain,
are commuting. This property can reduce the required number of steps within the
integrator significantly. If we consider two successive advection steps the
integrator is given by
\begin{align}
\begin{pmatrix}
    \vec X(t_0 + 2h)\\
    \vec V(t_0 + 2h)
\end{pmatrix} + O(h^2) 
                    =&\Bigl( \varphi_{t_0 + 2h,t_0 + 3h/2}^{\supidx{v_x}}\circ\varphi_{t_0 + 2h,t_0 + 3h/2}^{\supidx{v_y}}\circ\varphi_{t_0 + 2h,t_0 + 3h/2}^{\supidx{v_z}}\circ\nonumber\\
                    &\varphi_{t_0 + 2h,t_0+h}^{\supidx{z}}\circ\varphi_{t_0 + 2h,t_0+h}^{\supidx{y}}\circ\varphi_{t_0 + 2h,t_0+h}^{\supidx{x}}\circ\nonumber\\
                    &\varphi_{t_0 + 3h/2,t_0 + h/2}^{\supidx{v_x}}\circ\varphi_{t_0 + 3h/2,t_0 + h/2}^{\supidx{v_y}}\circ\varphi_{t_0 + 3h/2,t_0 + h/2}^{\supidx{v_z}}\circ\nonumber\\
                    &\varphi_{t_0 + h,t_0}^{\supidx{z}}\circ\varphi_{t_0 + h,t_0}^{\supidx{y}}\circ\varphi_{t_0 + h,t_0}^{\supidx{x}}\circ\nonumber\\
                    &\varphi_{t_0 + h/2,t_0}^{\supidx{v_z}}\circ\varphi_{t_0 + h/2,t_0}^{\supidx{v_y}}\circ\varphi_{t_0 + h/2,t_0}^{\supidx{v_x}}\Bigl)(\vecxv).
\end{align}
The above integrator has merged two half-time steps of
$\varphi^{\supidx{\vecv}}_{t0+h,t_0}$ using two properties. The first property
is the commutative property of two flows acting on different axes. The second
property is that we can add up to successive flows working on the same axis if
the intervals are adjacent to each other, often referred to as
"first-same-as-last" property. 
Merging these steps removes 30\% of the required operations during the advection
which is an important performance improvement, since advancing the distribution
function is the most expensive steps in solving the Vlasov equation. This was
discussed in detail in Schild et.al. \cite{Schild2024}.

Finally, we can reuse \cref{num:alg:SemiLagrangianPhysicalDomain} to actually
implement a semi-Lagrangian method on a rotating grid. The algorithm does not
change. The coefficients remain constant but depend on time which we have to
consider while solving \cref{num:equ:BacktracingTrajectory}
\begin{align}
\label{num:equ:BacktracingRotatingGridX}
X(0) &= x - \left(\int_{t_0 + h}^{t_0} \vec D_{\omega_c}^{-1}(t) \mathrm{d} t \vecv\right)_x\quad
Y(0) = y - \left(\int_{t_0 + h}^{t_0} \vec D_{\omega_c}^{-1}(t) \mathrm{d} t \vecv\right)_y\\
Z(0) &= z - \left(\int_{t_0 + h}^{t_0} \vec D_{\omega_c}^{-1}(t) \mathrm{d} t \vecv\right)_z\nonumber\\
\label{num:equ:BacktracingRotatingGridVx}
V_x(0) &= v_{x} - \left(\int_{t_0 + h}^{t_0} \vec D_{\omega_c}(t) \mathrm{d} t \vec E(\vecx)\right)_{v_x}\quad
V_y(0) = v_{y} - \left(\int_{t_0 + h}^{t_0} \vec D_{\omega_c}(t) \mathrm{d} t \vec E(\vecx)\right)_{v_y}\\
V_z(0) &= v_{z} - \left(\int_{t_0 + h}^{t_0} \vec D_{\omega_c}(t) \mathrm{d} t \vec E(\vecx)\right)_{v_z}\nonumber.
\end{align}

\section{Numerical comparison of the semi-Lagrangian method with and without a rotating grid}
\label{sec:NumericalResultsWAndWoRotatingGrid}
In this section we investigate the behavior of the rotating grid based on
different use cases. All simulations have been conducted using the performance 
portable BSL6D code~\cite{BSL6D} which is an open source project of the
numerical division of the Max-Plank-Institute of Plasma Physics.

\subsection{Solving the $\vecv \times \vec B_0$ term}
\label{subsec:NumericalExampleVlasovEquationVxB}
We start our investigation solving only the rotational part of the Vlasov
equation which is
\cref{num:equ:VlasovEquationRotationPart,num:equ:VlasovEquationRotationPartCompDomain}
and a basic proof of concept.  Using $\vec B_0=(0,0,1)$ we can solve these two
equation on a 2-D domain defined by $v_x$ and $v_y$. The low dimensionality also
allows us to fully visualize the distribution function which helps to understand
the behavior of the rotating grid.  We compare the results of a Strang-Splitting
with and without rotating velocity domain with the analytical solution. In the
rotating frame we first transform the solution back to the physical domain
and afterwards compare against the analytical solution. The transformation from
the computational domain into the physical domain based on the inverse mapping
of \cref{num:equ:PhysicalDomainToRotatingFrame}.

The initial condition for our test is given by
\begin{align}
\label{num:equ:NumericalAnalysisInitialConditionVelocity}
f_0(\vecv) = \dfrac{1}{\sqrt{2\pi}} \exp\left(-\dfrac{(\vecv-(1,0,0))^2}{2}\right),
\end{align}
where all occuring physical quantities have been normalized ($q=m=1$) to one. The
solution to the characteristic equations are given by a harmonic oscillator as
shown by Chen~\cite[Subsec.~2.2.1]{Chen2016}. We can use the analytical trajectories
to solve
\cref{num:equ:VlasovEquationRotationPart,num:equ:VlasovEquationRotationPartCompDomain}
using the semi-Lagrangian method and trace the grid points back to the initial
condition. With these trajectories the time dependent distribution function is given
by
\begin{align}
\label{num:equ:VlasovEquationRotationPartAnalyticalSolutionStrang}
f(\vecv,t) &= f_0(\vec D_{\omega_c}(t)\vecv)\\
\label{num:equ:VlasovEquationRotationPartAnalyticalSolutionRotatingGrid}
f(\tilde{\vecv},t) &= f(\vec D_{\omega_c}^{-1}(t) \tilde{\vecv},t) = f_0(\vec D_{\omega_c}(t)\vec D_{\omega_c}^{-1}(t)\tilde{\vecv}) = f_0(\tilde{\vecv})
\end{align}
where \cref{num:equ:VlasovEquationRotationPartAnalyticalSolutionStrang} provides
the solution in the physical domain while
\cref{num:equ:VlasovEquationRotationPartAnalyticalSolutionRotatingGrid} gives
the solution in the rotational domain which has been transformed into the
physical domain in the first equality. 

\begin{figure}
\resizebox{1\textwidth}{!}{
\input{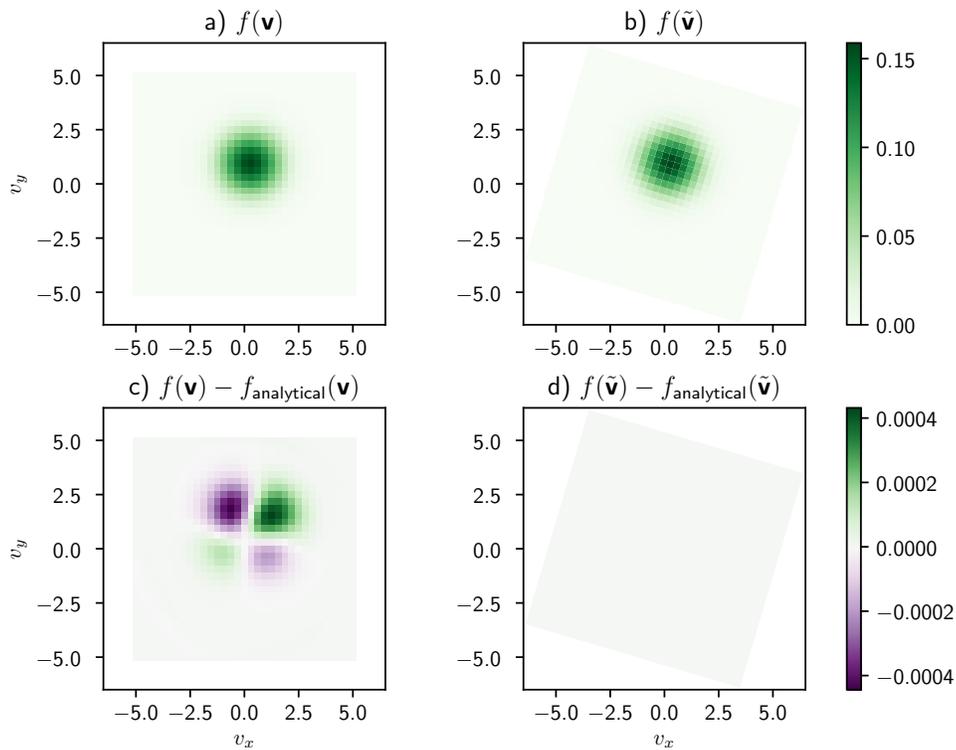}
}
\caption{\label{num:fig:rotatingGridComparison}Visualization of the solution of \cref{num:equ:VlasovEquationRotationPart} 
with a classical Strang-Splitting and on a rotating grid in the physical domain at $t=5/\omega_c$. The figures a) and b)
give the full solution while c) and d) show the difference of the simulation results and the analytical solution in
\cref{num:equ:VlasovEquationRotationPartAnalyticalSolutionStrang,num:equ:VlasovEquationRotationPartAnalyticalSolutionRotatingGrid}}
\end{figure}

The solution in the rotating frame is, as expected, a stationary solution. All
motion induced by the $\vecv\times \vec B_0$ term was removed from the Vlasov
equation. A visualization of the rotating state compared to the fixed state is
given in \cref{num:fig:rotatingGridComparison}. The difference of the analytical
and the simulation result is given in the L2 norm in
\cref{num:fig:rotatingGridInftyErr}. We do not plot the error of the rotating
grid since the solution only contains unit operations which do not change the 
initial condition. The error of the Strang-Splitting approach increases linearly
over time. The linear increase is superimposed by an oscillating component which
has its minima at symmetry position of the initial condition. At $t=(2n-1)\pi$
the solution is mirrored on the $v_y$ axis and at $t=2n\pi$ the analytical
solution is equivalent to the initial condition with $n\in\mathrm{N}$. At these
two time steps the oscillatory component of the error of the Strang-Splitting
approach is minimal with respect to the analytical solution.

One advantage of the consideration of exactly solvable problems is that these
problems provide perfect test cases for unit test in software applications. The
described setup of this subsection provides one example of a unit test which
continuously monitors the behavior of the BSL6D code\cite{BSL6D}.

\subsection{Solution of the Vlasov equation with constant background fields}
\label{subsec:NumericalExampleVlasovEquationConstBgFields}
In this example we focus on the convergence behavior of the integrators which
have been introduced in
\cref{subsec:IntegratorsForTheVlasovEquation,subsec:IntegratorsForTheVlasovEquationRotatingFrame}
for \cref{num:equ:VlasovEquation,num:equ:VlasovEquationRotatingFrame}.  We
extend the use case of the last subsection to the full Vlasov equation with
constant background fields using $\vec E(\vecx,t) = \vec E_0 = (E_0,0,0)$ with
$E_0=0.1$ and adding the advection part in the spatial domain. The electric
field is rather small which is consistent with our example in the next
subsection. The initial condition for the velocity space is again given by
\cref{num:equ:NumericalAnalysisInitialConditionVelocity}. The spatial domain is
initialized with a plane wave perturbation
\begin{align}
\label{num:equ:NumericalAnalysisInitialConditionSpatial}
f_0(\vecx) = 1 + \epsilon \sin(k_{0,x} x + k_{0,y} y)
\end{align}
using a small perturbation amplitude $\epsilon=0.1$ and the smallest modes
$k_{0,i}=2\pi/L_i$ with $i=x,y$ which can be represented on the spatial domain.
The initial condition is given by the product
$f_0(\vecxv)=f_0(\vecx)f_0(\vecv)$.

We can again utilize the semi-Lagrangian method to calculate the analytical
solution $f(\vecxv,t)$ by solving the characteristic equations. The trajectories
are solved e.g. by Chen~\cite[Subsec.~2.2.2]{Chen2016}. The solution in the rotating
frame is compared by first transforming the initial condition back to the
computational domain and afterwards tracing the grid points back to the initial
condition.  The time dependent distribution function is then given by
\begin{align}
f(&\vecxv,t) = f_0\left(\vecx + \int_{t}^{0} \vec D_{\omega_c}(t')\left(\vecv + \vec E_0\right)\mathrm{d}t' - \vec E_0 t ,\vec D_{\omega_c}(t)(\vecv + \vec E_0) - \vec E_0\right)\\
f(&\vecx,\tilde{\vecv},t) = f(\vecx,\vec D_{\omega_c}^{-1}(t)\tilde{\vecv},t)\nonumber\\
&= f_0\left(\vecx + \int_{t}^{0} \vec D_{\omega_c}(t')\left(\vec D_{\omega_c}^{-1}(t)\tilde{\vecv} + \vec E_0\right)\mathrm{d}t' - \vec E_0 t ,\vec D_{\omega_c}(t)(\vec D_{\omega_c}^{-1}(t)\tilde{\vecv} + \vec E_0) - \vec E_0\right),
\end{align}
where the first result is the solution in the physical domain while the second
result gives the solution in the rotational domain which has been transformed
into the physical domain in the first equivalence relation. 

\begin{figure}
\begin{center}
\resizebox{1\textwidth}{!}{
\def\mathdefault#1{#1}
\everymath=\expandafter{\the\everymath\displaystyle}
\input{Graphics/RotatingGrid_ExB_l2RelError.pgf}
}
\caption{\label{num:fig:rotatingGridInftyErr}Difference of the analytical
solution to simulation results for the test cases in
\cref{subsec:NumericalExampleVlasovEquationVxB,subsec:NumericalExampleVlasovEquationConstBgFields}
based on the L2 norm using $h=0.01/\omega_c$. The error is normalized with
respect to the perturbation of the analytical solution $\delta
f_\text{analytical}$. As is shown in \cref{num:fig:rotatingGridComparison} the
rotating grid has no numerical error for the $\vecv \times \vec B_0$ simulation,
which is therefore omitted in the plot.}
\end{center}
\end{figure}
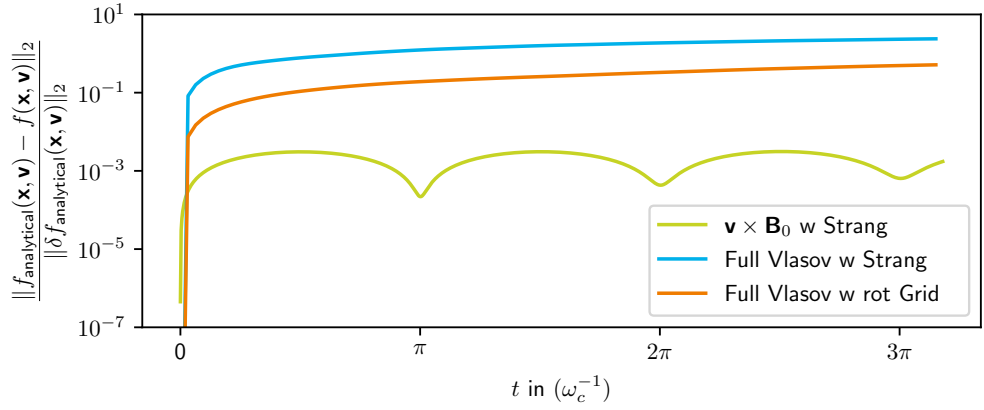

The difference of the analytical solution and the simulation result is also
given \cref{num:fig:rotatingGridInftyErr} using again the L2 norm. The
normalization is chosen such that we plot the relative error of the perturbation
$\delta f = f - 1$. Now also the error of the solution on the rotating grid
increases linearly.  But compared to the classical Strang-Splitting approach the
error is a magnitude smaller such that we can state that the rotating grid is
numerically advantageous compared to a pure Strang-Splitting approach.

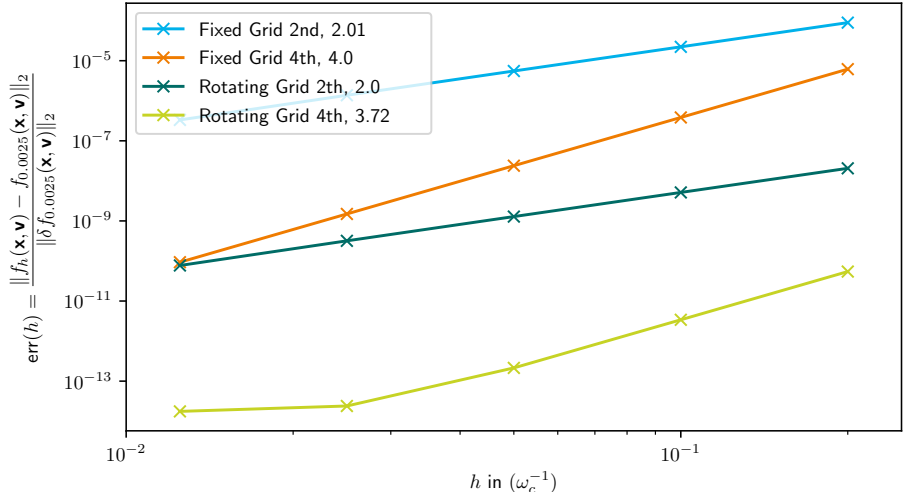
\begin{figure}
\begin{center}
\resizebox{1\textwidth}{!}{
\def\mathdefault#1{#1}
\everymath=\expandafter{\the\everymath\displaystyle}
\input{Graphics/RotatingGrid_convergenceRates.pgf}
}
\caption{\label{num:fig:rotatingGridConvergenceRates}Convergence rates for
integrators of
\cref{subsec:IntegratorsForTheVlasovEquation,subsec:IntegratorsForTheVlasovEquationRotatingFrame,subsec:NumericalExampleVlasovEquationConstBgFields}.
The convergence rate in the legend always omit the first data point of the
measurement. The errors are estimated against a solution that has been obtained
using a significantly smaller time step $h=0.0025/\omega_c$ and is referred to
as a converged solution. The comparison carried out at $t=9.0/\omega_c$. The
difference is normalized on the perturbation $\delta f$ of the converged
solution.}
\end{center}
\end{figure}

Finally, we validate the convergence rates which have been derived in
\cref{subsec:IntegratorsForTheVlasovEquation,subsec:IntegratorsForTheVlasovEquationRotatingFrame}.
These simulations have been based on trigonometric interpolation to allow for
larger time steps, which would not have been possible with the Lagrange
interpolation which does not allow such large time steps in the BSL6D
Code\cite{BSL6D}. The measured convergence rates are shown in
\cref{num:fig:rotatingGridConvergenceRates}. Additionally, to the convergence
rates of the Strang-Splitting approaches we added a fourth order splitting
schemes which can be constructed based on a Strang-Splitting taken from
Kraus et.al. ~\cite[p.~31]{Kraus2017}
\begin{align}
\varphi_{h,S4} = \varphi_{\gamma_1 h,S}\circ\varphi_{\gamma_2 h,S}\circ\varphi_{\gamma_1 h,S}
\end{align}
with
\begin{align}
\gamma_1 = \dfrac{1}{2-2^{1/3}} \qquad \gamma_2 = -\dfrac{2^{1/3}}{2-2^{1/3}}
\end{align}
where $\varphi_{\gamma_i h,S}$ is either the integrator given in
\cref{subsec:IntegratorsForTheVlasovEquation} or
\cref{subsec:IntegratorsForTheVlasovEquationRotatingFrame}. The measured
convergence rates match very well the expected convergence rates and are calculated 
using
\begin{align}
m = \dfrac{\log(\text{err}(0.2))-\log(\text{err}(0.025))}{\log(0.2)-\log(0.025)}.
\end{align}
Only the fourth order integrator combined with the rotating grid shows
deviations from the expected convergence rates for small $h$. Since the
difference between the converged solution and the simulation result is rather
small with a difference of $10^{-13}$ the deviation can be justified by other
discretization, rounding, or interpolation errors which dominate in this error
regime.

As in the last subsection also this setup provides us with a perfect unit test
which is used to continuously validate the behavior of our the implementation in 
\cite{BSL6D}.

\subsection{Coupling the Vlasov equation to the quasi-neutrality equation}
\label{subsec:NonlinearVlasovExamples}
\paragraph{Stable neutralized ion Bernstein waves: Dispersion relation} In this
last subsection we consider nonlinear examples described by
\eqref{num:equ:QuasiNeutralWithAdiabaticElectrons}. We normalized physical
quantities ($e=T=m=1$) in these equations. The electric field is coupled to the
distribution function through the quasi-neutrality condition with adiabatic
electrons in \cref{num:equ:QuasiNeutralWithAdiabaticElectrons}.

In the first example we reproduce the dispersion relation of neutralizing ion
Bernstein waves (nIBW) \cite{brambilla_kinetic_1998} which have been one central aspects of
the study of the limits of gyrokinetics in \cite{RaethThesis2023}.  The example can be solved as a
3-D problem which consists of the dimensions $y,v_x,v_y$.  The velocity space
contains the full rotation and the dispersion relation is reproduced as
$\omega(k_y)$. We choose the initial condition to specifically excite nIBWs in
the $y$ dimension of our simulation
\begin{align}
f_0(\vecxv) = f_0(\vecv) \biggl[ 1 + \alpha \sum_{m_{k_y}=1}^{m_{\text{max}}} \sum_{p=0}^{p_{\text{max}}} &J_p(k_y v_\perp)\cdot \text{min}\left(\dfrac{1}{e^{-k_y^2} I_p(k_y^2)},0.01(p+1)^{1/3}\right)\cdot \nonumber\\
                                                                                                            &\text{Re}(e^{\text{i}v_\perp k_y \sin(\gamma) - p\gamma + k_y y}) \biggr]
\end{align}
Here $ I_p(\cdot)$ and $J_p(\cdot)$ are the modified cylindrical Bessel
functions and Bessel Functions of the first kind, respectively. Also, the
perpendicular velocity $v_\perp^2 = v_x^2 + v_y^2$ and the angle
$\gamma\sphericalangle (v_\perp,k_y)$ are needed for the initialization.
Finally, $\alpha$, $m_\text{max}$, and $p_\text{max}$ are the perturbation
amplitude, the maximal mode and the maximum order of Bessel functions,
respectively. The initialization is based on the analytical solution of nIBWs
for this numerical example.

\begin{figure}
\begin{center}
\includegraphics[scale=0.9]{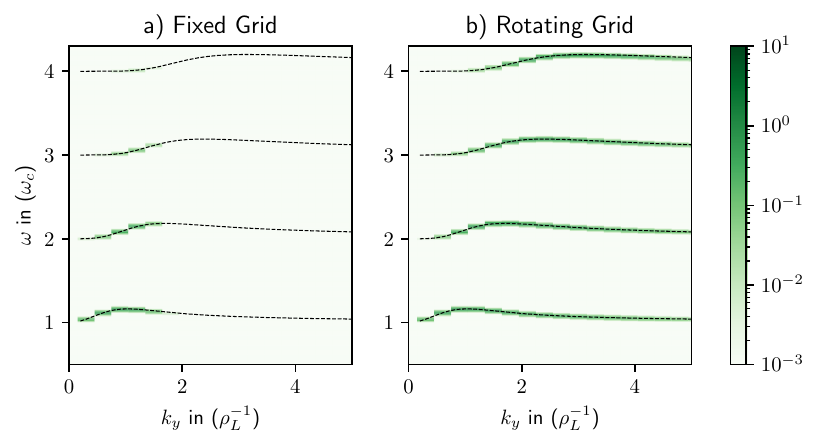}
\caption{\label{num:fig:compareRotFixedGridNIBW}Dispersion relation of
neutralizing ion Bernstein waves with a) fixed velocity grid and b) rotating
velocity grid. The black dashed lines represent the analytical solution of the
example. The length scale is normalized to the Larmor radius $\rho_L$.}
\end{center}
\end{figure}

After the initialization the simulation is executed to $t=1000$ using
$h=0.05$. A Fourier transform is applied to the resulting particle density in
space and time $n(\vecx,t)$ to obtain the dispersion relation which is plotted
in \cref{num:fig:compareRotFixedGridNIBW}.  We can observe a clear quantitative
and qualitative difference in our results.  The branches of the dispersion
relation with a classical Strang splitting are only visible in the range
$k_y\in(0,\rho_L^{-1})$ for the first two harmonics of the gyrofrequency
$\omega_c$.  With the rotating grid we can reproduce the branches of the
dispersion relation within the full domain which has been plotted. For higher
waves the Lagrange interpolation has damping effects which removes all
perturbations of the distribution function. This is not visible in this domain
but was shown in \cite[p.11-12]{Schild2024}.

Our conclusion based on this nonlinear numerical example is that the rotating
grid is clearly advantageous compared to simple splitting approaches and
furthermore allows for a significant reduction of splitting steps due to merging
of splitted steps as discussed in
\cref{subsec:IntegratorsForTheVlasovEquationRotatingFrame}.

\begin{figure}
    \begin{center}
    \resizebox{1\textwidth}{!}{
\def\mathdefault#1{#1}
\everymath=\expandafter{\the\everymath\displaystyle}
    \input{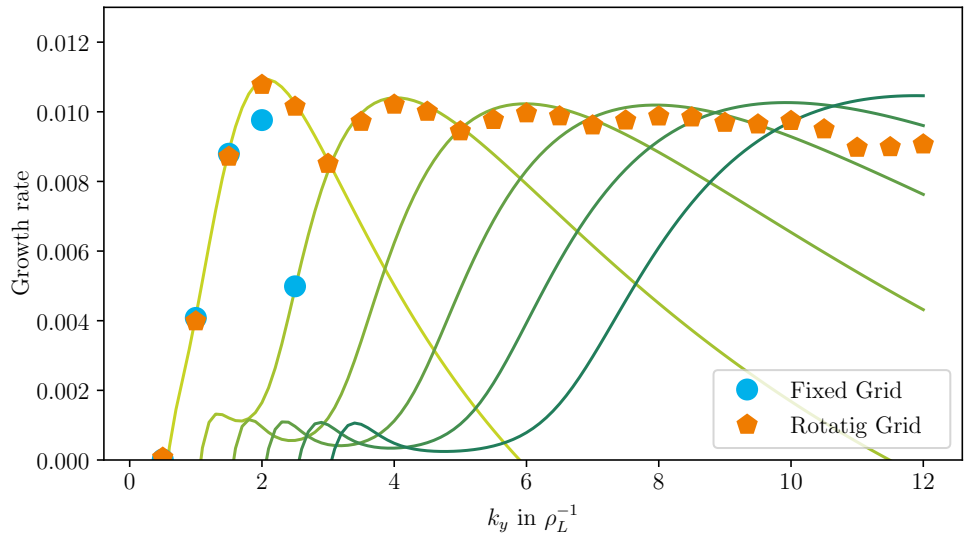}}
    \caption{\label{num:fig:compareRotFixedGridUnstableNIBW} Growth rate of 
    neutralizing ion Bernstein waves with a) fixed velocity grid and Strang splitting and b) rotating
    velocity grid. The solid lines represent the analytical solution of the growth rates for the first six harmonics of the IBWs with increasing frequency from light to dark green.}
    \end{center}
    \end{figure}

\paragraph{Unstable neutralized ion Bernstein waves: Growth rate}
The unstable neutralized ion Bernstein waves are also based on
\cref{num:equ:QuasiNeutralWithAdiabaticElectrons}. But the simulation is
executed in a different setup. Neutralized ion Bernstein waves can be
destabilized through the imposition of density and temperature gradients
\cite{RaethPhysRevLett24}. To achieve this, we introduce a right-hand side term
to \cref{num:equ:VlasovEquation}, incorporating temperature and density
gradients, as outlined in \cite{RaethPhysPlasmas24}
\begin{align}
	\partial_t f + \vec v\cdot \nablax f + \left[-\nablax \phi+(\vec v\times \vec B_0)\right] \cdot \nablav f =  \vecv^*\cdot \nablax\phi f_M,
\end{align}
where $f_M$ a constant background distribution function $f = f_M + \delta f$
which introduces a density and temperature gradient on the background distribution.
The source on the right-hand side is given by 
\begin{align}
	\vecv^* = \vec B_0 \times \frac{ \nabla n}{n} - \vec B_0 \times \left( \frac{ \nabla T}{T} \frac{3 -  (\perpc v ^2 + \parac v^2)}{2 }\right).
\end{align}
Here we use $v_\vert^2=v_x^2+v_y^2$ which is perpendicular to the magnetic field. 
The parameters for the gradients are
\begin{align}
	\kappa_n = \frac{\partial_x n}{n} = 0.44 ; \hspace{2cm} \kappa_T = \frac{\partial_x T}{T} = 0.36. \label{num:equ:gradient_parameter}
\end{align}
A resolution in configuration space with $N = 1\times 256\times 8\times 33\times
33\times 33$ has been chosen for a box with length $L = \pi\times 4\pi\times
80\pi$. The simulation has been performed with a time step of $\Delta t  =
0.005$.  Figure \ref{num:fig:compareRotFixedGridUnstableNIBW} illustrates a
comparison of growth rates between simulations utilizing the rotating grid and
the Strang splitting. Additionally, the analytical dispersion relation from
\cite{RaethPhysPlasmas24} for parameter given in
\eqref{num:equ:gradient_parameter} is included. Although the simulation with the
rotating grid slightly deviates from the analytical results due to numerical
damping from spatial advection interpolations, it accurately reproduces the
correct growth rate for a significant range of wave numbers. In contrast, in the
simulation employing Strang splitting, only the growth rates for the first two
wave numbers are accurately reproduced.

\section{Acknowledgments} 
Computations have been performed on the HPC system Raven at the Max Planck
Computing and Data Facility. Additionally, we thank Omar Maj and Tileuzhan
Mukhamet for fruitful discussions on the convergence analysis and coordinate
transformation.

This work has been carried out within the framework of the EUROfusion
Consortium, funded by the European Union via the Euratom Research and Training
Programme (Grant Agreement No 101052200 — EUROfusion). Views and opinions
expressed are however those of the author(s) only and do not necessarily reflect
those of the European Union or the European Commission. Neither the European
Union nor the European Commission can be held responsible for them.

\appendix
\section{Interpolation Methods}
\label{app:InterpolationMethods}
In this paper we use two different interpolations which we briefly introduce.
Both assume equidistant grid points $x_j=j\Delta x$, $j=1,\cdots,N$ and
interpolate a 1-D function.

The first interpolation is the Lagrange interpolation $L(x)$ which has
performance advantages due to its locality as discussed in \cite{Schild2024}.
The locality has the drawback that the interpolated point has to be centered by
the Lagrange interpolation stencil.  This can introduce an implementation based
CFL condition. Therefore, we assume that the interpolation shift $\alpha$ is
smaller than the spacing of the grid $\Delta x$.  Furthermore, we denote by
$l^q_i$ the Lagrange-polynomials of order $(q-1)$ with $q$ nodes in the
interpolant
\begin{itemize}
    \item For an odd number q, the interpolant is given by
    \begin{align}
        L(x_j + \alpha) = \sum_{i=j-(q-1)/2}^{j+(q-1)/2}l^q_i(\alpha)f(x_i)
    \end{align}
    \item For an even number q, the interpolation stencil is centered around
    the interpolated point $x_j + \alpha$, such that the interpolation is given
    by
    \begin{align}
        L(x_j + \alpha) = \begin{cases}
            \sum_{i=j-q/2}^{j+q/2-1} l^q_i(\alpha)f(x_i)\\
            \sum_{i=j-q/2+1}^{j+q/2} l^q_i(\alpha)f(x_i)
        \end{cases}.
    \end{align}
\end{itemize}

The second interpolation method is a trigonometric interpolation $T(x)$ which is a
global formula but more accurate compared to the Lagrange interpolation. The
implemented interpolant is given by
\begin{itemize}
    \item For an even number $N$ of grid points
      \begin{align}
          T(x) = \sum_{j=1}^{N}\dfrac{\text{sinc}(1/2 N((x - x_j)))}{\text{sinc}(1/2(x-x_j))}\cos(1/2(x-x_j))f(x_j)
      \end{align}
    \item For an odd number $N$ of grid points
      \begin{align}
          T(x) = \sum_{j=1}^{N}\dfrac{\text{sinc}(1/2 N((x - x_j)))}{\text{sinc}(1/2(x-x_j))}f(x_j)
      \end{align}
\end{itemize}
The trigonometric interpolation could also be based on a Fast Fourier Transform
(FFT). 

\bibliographystyle{siamplain}
\bibliography{Bibliographie}

\end{document}

%% file: Graphics/RotatingGrid_ExB_l2RelError.pgf
\begingroup%
\makeatletter%
\begin{pgfpicture}%
\pgfpathrectangle{\pgfpointorigin}{\pgfqpoint{5.905512in}{2.755906in}}%
\pgfusepath{use as bounding box, clip}%
\begin{pgfscope}%
\pgfsetbuttcap%
\pgfsetmiterjoin%
\definecolor{currentfill}{rgb}{1.000000,1.000000,1.000000}%
\pgfsetfillcolor{currentfill}%
\pgfsetlinewidth{0.000000pt}%
\definecolor{currentstroke}{rgb}{1.000000,1.000000,1.000000}%
\pgfsetstrokecolor{currentstroke}%
\pgfsetdash{}{0pt}%
\pgfpathmoveto{\pgfqpoint{0.000000in}{0.000000in}}%
\pgfpathlineto{\pgfqpoint{5.905512in}{0.000000in}}%
\pgfpathlineto{\pgfqpoint{5.905512in}{2.755906in}}%
\pgfpathlineto{\pgfqpoint{0.000000in}{2.755906in}}%
\pgfpathlineto{\pgfqpoint{0.000000in}{0.000000in}}%
\pgfpathclose%
\pgfusepath{fill}%
\end{pgfscope}%
\begin{pgfscope}%
\pgfsetbuttcap%
\pgfsetmiterjoin%
\definecolor{currentfill}{rgb}{1.000000,1.000000,1.000000}%
\pgfsetfillcolor{currentfill}%
\pgfsetlinewidth{0.000000pt}%
\definecolor{currentstroke}{rgb}{0.000000,0.000000,0.000000}%
\pgfsetstrokecolor{currentstroke}%
\pgfsetstrokeopacity{0.000000}%
\pgfsetdash{}{0pt}%
\pgfpathmoveto{\pgfqpoint{0.776331in}{0.668778in}}%
\pgfpathlineto{\pgfqpoint{5.755512in}{0.668778in}}%
\pgfpathlineto{\pgfqpoint{5.755512in}{2.527906in}}%
\pgfpathlineto{\pgfqpoint{0.776331in}{2.527906in}}%
\pgfpathlineto{\pgfqpoint{0.776331in}{0.668778in}}%
\pgfpathclose%
\pgfusepath{fill}%
\end{pgfscope}%
\begin{pgfscope}%
\pgfsetbuttcap%
\pgfsetroundjoin%
\definecolor{currentfill}{rgb}{0.000000,0.000000,0.000000}%
\pgfsetfillcolor{currentfill}%
\pgfsetlinewidth{0.803000pt}%
\definecolor{currentstroke}{rgb}{0.000000,0.000000,0.000000}%
\pgfsetstrokecolor{currentstroke}%
\pgfsetdash{}{0pt}%
\pgfsys@defobject{currentmarker}{\pgfqpoint{0.000000in}{-0.048611in}}{\pgfqpoint{0.000000in}{0.000000in}}{%
\pgfpathmoveto{\pgfqpoint{0.000000in}{0.000000in}}%
\pgfpathlineto{\pgfqpoint{0.000000in}{-0.048611in}}%
\pgfusepath{stroke,fill}%
}%
\begin{pgfscope}%
\pgfsys@transformshift{1.002657in}{0.668778in}%
\pgfsys@useobject{currentmarker}{}%
\end{pgfscope}%
\end{pgfscope}%
\begin{pgfscope}%
\definecolor{textcolor}{rgb}{0.000000,0.000000,0.000000}%
\pgfsetstrokecolor{textcolor}%
\pgfsetfillcolor{textcolor}%
\pgftext[x=1.002657in,y=0.571556in,,top]{\color{textcolor}{\sffamily\fontsize{10.000000}{12.000000}\selectfont\catcode`\^=\active\def^{\ifmmode\sp\else\^{}\fi}\catcode`\%=\active\def
\end{pgfscope}%
\begin{pgfscope}%
\pgfsetbuttcap%
\pgfsetroundjoin%
\definecolor{currentfill}{rgb}{0.000000,0.000000,0.000000}%
\pgfsetfillcolor{currentfill}%
\pgfsetlinewidth{0.803000pt}%
\definecolor{currentstroke}{rgb}{0.000000,0.000000,0.000000}%
\pgfsetstrokecolor{currentstroke}%
\pgfsetdash{}{0pt}%
\pgfsys@defobject{currentmarker}{\pgfqpoint{0.000000in}{-0.048611in}}{\pgfqpoint{0.000000in}{0.000000in}}{%
\pgfpathmoveto{\pgfqpoint{0.000000in}{0.000000in}}%
\pgfpathlineto{\pgfqpoint{0.000000in}{-0.048611in}}%
\pgfusepath{stroke,fill}%
}%
\begin{pgfscope}%
\pgfsys@transformshift{2.426131in}{0.668778in}%
\pgfsys@useobject{currentmarker}{}%
\end{pgfscope}%
\end{pgfscope}%
\begin{pgfscope}%
\definecolor{textcolor}{rgb}{0.000000,0.000000,0.000000}%
\pgfsetstrokecolor{textcolor}%
\pgfsetfillcolor{textcolor}%
\pgftext[x=2.426131in,y=0.571556in,,top]{\color{textcolor}{\sffamily\fontsize{10.000000}{12.000000}\selectfont\catcode`\^=\active\def^{\ifmmode\sp\else\^{}\fi}\catcode`\%=\active\def
\end{pgfscope}%
\begin{pgfscope}%
\pgfsetbuttcap%
\pgfsetroundjoin%
\definecolor{currentfill}{rgb}{0.000000,0.000000,0.000000}%
\pgfsetfillcolor{currentfill}%
\pgfsetlinewidth{0.803000pt}%
\definecolor{currentstroke}{rgb}{0.000000,0.000000,0.000000}%
\pgfsetstrokecolor{currentstroke}%
\pgfsetdash{}{0pt}%
\pgfsys@defobject{currentmarker}{\pgfqpoint{0.000000in}{-0.048611in}}{\pgfqpoint{0.000000in}{0.000000in}}{%
\pgfpathmoveto{\pgfqpoint{0.000000in}{0.000000in}}%
\pgfpathlineto{\pgfqpoint{0.000000in}{-0.048611in}}%
\pgfusepath{stroke,fill}%
}%
\begin{pgfscope}%
\pgfsys@transformshift{3.849606in}{0.668778in}%
\pgfsys@useobject{currentmarker}{}%
\end{pgfscope}%
\end{pgfscope}%
\begin{pgfscope}%
\definecolor{textcolor}{rgb}{0.000000,0.000000,0.000000}%
\pgfsetstrokecolor{textcolor}%
\pgfsetfillcolor{textcolor}%
\pgftext[x=3.849606in,y=0.571556in,,top]{\color{textcolor}{\sffamily\fontsize{10.000000}{12.000000}\selectfont\catcode`\^=\active\def^{\ifmmode\sp\else\^{}\fi}\catcode`\%=\active\def
\end{pgfscope}%
\begin{pgfscope}%
\pgfsetbuttcap%
\pgfsetroundjoin%
\definecolor{currentfill}{rgb}{0.000000,0.000000,0.000000}%
\pgfsetfillcolor{currentfill}%
\pgfsetlinewidth{0.803000pt}%
\definecolor{currentstroke}{rgb}{0.000000,0.000000,0.000000}%
\pgfsetstrokecolor{currentstroke}%
\pgfsetdash{}{0pt}%
\pgfsys@defobject{currentmarker}{\pgfqpoint{0.000000in}{-0.048611in}}{\pgfqpoint{0.000000in}{0.000000in}}{%
\pgfpathmoveto{\pgfqpoint{0.000000in}{0.000000in}}%
\pgfpathlineto{\pgfqpoint{0.000000in}{-0.048611in}}%
\pgfusepath{stroke,fill}%
}%
\begin{pgfscope}%
\pgfsys@transformshift{5.273080in}{0.668778in}%
\pgfsys@useobject{currentmarker}{}%
\end{pgfscope}%
\end{pgfscope}%
\begin{pgfscope}%
\definecolor{textcolor}{rgb}{0.000000,0.000000,0.000000}%
\pgfsetstrokecolor{textcolor}%
\pgfsetfillcolor{textcolor}%
\pgftext[x=5.273080in,y=0.571556in,,top]{\color{textcolor}{\sffamily\fontsize{10.000000}{12.000000}\selectfont\catcode`\^=\active\def^{\ifmmode\sp\else\^{}\fi}\catcode`\%=\active\def
\end{pgfscope}%
\begin{pgfscope}%
\definecolor{textcolor}{rgb}{0.000000,0.000000,0.000000}%
\pgfsetstrokecolor{textcolor}%
\pgfsetfillcolor{textcolor}%
\pgftext[x=3.265921in,y=0.381587in,,top]{\color{textcolor}{\sffamily\fontsize{10.000000}{12.000000}\selectfont\catcode`\^=\active\def^{\ifmmode\sp\else\^{}\fi}\catcode`\%=\active\def
\end{pgfscope}%
\begin{pgfscope}%
\pgfsetbuttcap%
\pgfsetroundjoin%
\definecolor{currentfill}{rgb}{0.000000,0.000000,0.000000}%
\pgfsetfillcolor{currentfill}%
\pgfsetlinewidth{0.803000pt}%
\definecolor{currentstroke}{rgb}{0.000000,0.000000,0.000000}%
\pgfsetstrokecolor{currentstroke}%
\pgfsetdash{}{0pt}%
\pgfsys@defobject{currentmarker}{\pgfqpoint{-0.048611in}{0.000000in}}{\pgfqpoint{-0.000000in}{0.000000in}}{%
\pgfpathmoveto{\pgfqpoint{-0.000000in}{0.000000in}}%
\pgfpathlineto{\pgfqpoint{-0.048611in}{0.000000in}}%
\pgfusepath{stroke,fill}%
}%
\begin{pgfscope}%
\pgfsys@transformshift{0.776331in}{0.668778in}%
\pgfsys@useobject{currentmarker}{}%
\end{pgfscope}%
\end{pgfscope}%
\begin{pgfscope}%
\definecolor{textcolor}{rgb}{0.000000,0.000000,0.000000}%
\pgfsetstrokecolor{textcolor}%
\pgfsetfillcolor{textcolor}%
\pgftext[x=0.391106in, y=0.616016in, left, base]{\color{textcolor}{\sffamily\fontsize{10.000000}{12.000000}\selectfont\catcode`\^=\active\def^{\ifmmode\sp\else\^{}\fi}\catcode`\%=\active\def
\end{pgfscope}%
\begin{pgfscope}%
\pgfsetbuttcap%
\pgfsetroundjoin%
\definecolor{currentfill}{rgb}{0.000000,0.000000,0.000000}%
\pgfsetfillcolor{currentfill}%
\pgfsetlinewidth{0.803000pt}%
\definecolor{currentstroke}{rgb}{0.000000,0.000000,0.000000}%
\pgfsetstrokecolor{currentstroke}%
\pgfsetdash{}{0pt}%
\pgfsys@defobject{currentmarker}{\pgfqpoint{-0.048611in}{0.000000in}}{\pgfqpoint{-0.000000in}{0.000000in}}{%
\pgfpathmoveto{\pgfqpoint{-0.000000in}{0.000000in}}%
\pgfpathlineto{\pgfqpoint{-0.048611in}{0.000000in}}%
\pgfusepath{stroke,fill}%
}%
\begin{pgfscope}%
\pgfsys@transformshift{0.776331in}{1.133560in}%
\pgfsys@useobject{currentmarker}{}%
\end{pgfscope}%
\end{pgfscope}%
\begin{pgfscope}%
\definecolor{textcolor}{rgb}{0.000000,0.000000,0.000000}%
\pgfsetstrokecolor{textcolor}%
\pgfsetfillcolor{textcolor}%
\pgftext[x=0.391106in, y=1.080798in, left, base]{\color{textcolor}{\sffamily\fontsize{10.000000}{12.000000}\selectfont\catcode`\^=\active\def^{\ifmmode\sp\else\^{}\fi}\catcode`\%=\active\def
\end{pgfscope}%
\begin{pgfscope}%
\pgfsetbuttcap%
\pgfsetroundjoin%
\definecolor{currentfill}{rgb}{0.000000,0.000000,0.000000}%
\pgfsetfillcolor{currentfill}%
\pgfsetlinewidth{0.803000pt}%
\definecolor{currentstroke}{rgb}{0.000000,0.000000,0.000000}%
\pgfsetstrokecolor{currentstroke}%
\pgfsetdash{}{0pt}%
\pgfsys@defobject{currentmarker}{\pgfqpoint{-0.048611in}{0.000000in}}{\pgfqpoint{-0.000000in}{0.000000in}}{%
\pgfpathmoveto{\pgfqpoint{-0.000000in}{0.000000in}}%
\pgfpathlineto{\pgfqpoint{-0.048611in}{0.000000in}}%
\pgfusepath{stroke,fill}%
}%
\begin{pgfscope}%
\pgfsys@transformshift{0.776331in}{1.598342in}%
\pgfsys@useobject{currentmarker}{}%
\end{pgfscope}%
\end{pgfscope}%
\begin{pgfscope}%
\definecolor{textcolor}{rgb}{0.000000,0.000000,0.000000}%
\pgfsetstrokecolor{textcolor}%
\pgfsetfillcolor{textcolor}%
\pgftext[x=0.391106in, y=1.545580in, left, base]{\color{textcolor}{\sffamily\fontsize{10.000000}{12.000000}\selectfont\catcode`\^=\active\def^{\ifmmode\sp\else\^{}\fi}\catcode`\%=\active\def
\end{pgfscope}%
\begin{pgfscope}%
\pgfsetbuttcap%
\pgfsetroundjoin%
\definecolor{currentfill}{rgb}{0.000000,0.000000,0.000000}%
\pgfsetfillcolor{currentfill}%
\pgfsetlinewidth{0.803000pt}%
\definecolor{currentstroke}{rgb}{0.000000,0.000000,0.000000}%
\pgfsetstrokecolor{currentstroke}%
\pgfsetdash{}{0pt}%
\pgfsys@defobject{currentmarker}{\pgfqpoint{-0.048611in}{0.000000in}}{\pgfqpoint{-0.000000in}{0.000000in}}{%
\pgfpathmoveto{\pgfqpoint{-0.000000in}{0.000000in}}%
\pgfpathlineto{\pgfqpoint{-0.048611in}{0.000000in}}%
\pgfusepath{stroke,fill}%
}%
\begin{pgfscope}%
\pgfsys@transformshift{0.776331in}{2.063124in}%
\pgfsys@useobject{currentmarker}{}%
\end{pgfscope}%
\end{pgfscope}%
\begin{pgfscope}%
\definecolor{textcolor}{rgb}{0.000000,0.000000,0.000000}%
\pgfsetstrokecolor{textcolor}%
\pgfsetfillcolor{textcolor}%
\pgftext[x=0.391106in, y=2.010362in, left, base]{\color{textcolor}{\sffamily\fontsize{10.000000}{12.000000}\selectfont\catcode`\^=\active\def^{\ifmmode\sp\else\^{}\fi}\catcode`\%=\active\def
\end{pgfscope}%
\begin{pgfscope}%
\pgfsetbuttcap%
\pgfsetroundjoin%
\definecolor{currentfill}{rgb}{0.000000,0.000000,0.000000}%
\pgfsetfillcolor{currentfill}%
\pgfsetlinewidth{0.803000pt}%
\definecolor{currentstroke}{rgb}{0.000000,0.000000,0.000000}%
\pgfsetstrokecolor{currentstroke}%
\pgfsetdash{}{0pt}%
\pgfsys@defobject{currentmarker}{\pgfqpoint{-0.048611in}{0.000000in}}{\pgfqpoint{-0.000000in}{0.000000in}}{%
\pgfpathmoveto{\pgfqpoint{-0.000000in}{0.000000in}}%
\pgfpathlineto{\pgfqpoint{-0.048611in}{0.000000in}}%
\pgfusepath{stroke,fill}%
}%
\begin{pgfscope}%
\pgfsys@transformshift{0.776331in}{2.527906in}%
\pgfsys@useobject{currentmarker}{}%
\end{pgfscope}%
\end{pgfscope}%
\begin{pgfscope}%
\definecolor{textcolor}{rgb}{0.000000,0.000000,0.000000}%
\pgfsetstrokecolor{textcolor}%
\pgfsetfillcolor{textcolor}%
\pgftext[x=0.477912in, y=2.475144in, left, base]{\color{textcolor}{\sffamily\fontsize{10.000000}{12.000000}\selectfont\catcode`\^=\active\def^{\ifmmode\sp\else\^{}\fi}\catcode`\%=\active\def
\end{pgfscope}%
\begin{pgfscope}%
\definecolor{textcolor}{rgb}{0.000000,0.000000,0.000000}%
\pgfsetstrokecolor{textcolor}%
\pgfsetfillcolor{textcolor}%
\pgftext[x=0.335550in,y=1.598342in,,bottom,rotate=90.000000]{\color{textcolor}{\sffamily\fontsize{10.000000}{12.000000}\selectfont\catcode`\^=\active\def^{\ifmmode\sp\else\^{}\fi}\catcode`\%=\active\def
\end{pgfscope}%
\begin{pgfscope}%
\pgfpathrectangle{\pgfqpoint{0.776331in}{0.668778in}}{\pgfqpoint{4.979181in}{1.859128in}}%
\pgfusepath{clip}%
\pgfsetrectcap%
\pgfsetroundjoin%
\pgfsetlinewidth{1.505625pt}%
\definecolor{currentstroke}{rgb}{0.776471,0.827451,0.145098}%
\pgfsetstrokecolor{currentstroke}%
\pgfsetdash{}{0pt}%
\pgfpathmoveto{\pgfqpoint{1.002657in}{0.820818in}}%
\pgfpathlineto{\pgfqpoint{1.007188in}{1.246366in}}%
\pgfpathlineto{\pgfqpoint{1.011719in}{1.316306in}}%
\pgfpathlineto{\pgfqpoint{1.016250in}{1.357216in}}%
\pgfpathlineto{\pgfqpoint{1.025312in}{1.408740in}}%
\pgfpathlineto{\pgfqpoint{1.034374in}{1.442655in}}%
\pgfpathlineto{\pgfqpoint{1.043437in}{1.467961in}}%
\pgfpathlineto{\pgfqpoint{1.052499in}{1.488143in}}%
\pgfpathlineto{\pgfqpoint{1.066092in}{1.512350in}}%
\pgfpathlineto{\pgfqpoint{1.079685in}{1.531783in}}%
\pgfpathlineto{\pgfqpoint{1.097809in}{1.552846in}}%
\pgfpathlineto{\pgfqpoint{1.115934in}{1.570125in}}%
\pgfpathlineto{\pgfqpoint{1.138589in}{1.588055in}}%
\pgfpathlineto{\pgfqpoint{1.161244in}{1.603059in}}%
\pgfpathlineto{\pgfqpoint{1.188430in}{1.618253in}}%
\pgfpathlineto{\pgfqpoint{1.220148in}{1.633105in}}%
\pgfpathlineto{\pgfqpoint{1.256396in}{1.647248in}}%
\pgfpathlineto{\pgfqpoint{1.297176in}{1.660428in}}%
\pgfpathlineto{\pgfqpoint{1.342486in}{1.672459in}}%
\pgfpathlineto{\pgfqpoint{1.392328in}{1.683187in}}%
\pgfpathlineto{\pgfqpoint{1.446701in}{1.692465in}}%
\pgfpathlineto{\pgfqpoint{1.505605in}{1.700140in}}%
\pgfpathlineto{\pgfqpoint{1.569039in}{1.706027in}}%
\pgfpathlineto{\pgfqpoint{1.632474in}{1.709713in}}%
\pgfpathlineto{\pgfqpoint{1.700440in}{1.711395in}}%
\pgfpathlineto{\pgfqpoint{1.768406in}{1.710807in}}%
\pgfpathlineto{\pgfqpoint{1.836372in}{1.707921in}}%
\pgfpathlineto{\pgfqpoint{1.899807in}{1.703034in}}%
\pgfpathlineto{\pgfqpoint{1.958711in}{1.696412in}}%
\pgfpathlineto{\pgfqpoint{2.017614in}{1.687498in}}%
\pgfpathlineto{\pgfqpoint{2.067456in}{1.677847in}}%
\pgfpathlineto{\pgfqpoint{2.112767in}{1.667040in}}%
\pgfpathlineto{\pgfqpoint{2.153546in}{1.655275in}}%
\pgfpathlineto{\pgfqpoint{2.189795in}{1.642786in}}%
\pgfpathlineto{\pgfqpoint{2.226043in}{1.627843in}}%
\pgfpathlineto{\pgfqpoint{2.257760in}{1.612130in}}%
\pgfpathlineto{\pgfqpoint{2.284947in}{1.596047in}}%
\pgfpathlineto{\pgfqpoint{2.307602in}{1.580189in}}%
\pgfpathlineto{\pgfqpoint{2.330257in}{1.561351in}}%
\pgfpathlineto{\pgfqpoint{2.348382in}{1.543452in}}%
\pgfpathlineto{\pgfqpoint{2.366506in}{1.522302in}}%
\pgfpathlineto{\pgfqpoint{2.384630in}{1.497320in}}%
\pgfpathlineto{\pgfqpoint{2.416348in}{1.451639in}}%
\pgfpathlineto{\pgfqpoint{2.420879in}{1.447770in}}%
\pgfpathlineto{\pgfqpoint{2.425410in}{1.445569in}}%
\pgfpathlineto{\pgfqpoint{2.429941in}{1.445271in}}%
\pgfpathlineto{\pgfqpoint{2.434472in}{1.446910in}}%
\pgfpathlineto{\pgfqpoint{2.439003in}{1.450307in}}%
\pgfpathlineto{\pgfqpoint{2.448065in}{1.460989in}}%
\pgfpathlineto{\pgfqpoint{2.497907in}{1.531739in}}%
\pgfpathlineto{\pgfqpoint{2.516031in}{1.551431in}}%
\pgfpathlineto{\pgfqpoint{2.534155in}{1.568178in}}%
\pgfpathlineto{\pgfqpoint{2.556810in}{1.585922in}}%
\pgfpathlineto{\pgfqpoint{2.579466in}{1.600959in}}%
\pgfpathlineto{\pgfqpoint{2.606652in}{1.616302in}}%
\pgfpathlineto{\pgfqpoint{2.638369in}{1.631372in}}%
\pgfpathlineto{\pgfqpoint{2.674618in}{1.645769in}}%
\pgfpathlineto{\pgfqpoint{2.715397in}{1.659215in}}%
\pgfpathlineto{\pgfqpoint{2.760708in}{1.671513in}}%
\pgfpathlineto{\pgfqpoint{2.810550in}{1.682503in}}%
\pgfpathlineto{\pgfqpoint{2.864922in}{1.692040in}}%
\pgfpathlineto{\pgfqpoint{2.923826in}{1.699973in}}%
\pgfpathlineto{\pgfqpoint{2.987261in}{1.706123in}}%
\pgfpathlineto{\pgfqpoint{3.050696in}{1.710066in}}%
\pgfpathlineto{\pgfqpoint{3.118662in}{1.712028in}}%
\pgfpathlineto{\pgfqpoint{3.186628in}{1.711738in}}%
\pgfpathlineto{\pgfqpoint{3.254594in}{1.709182in}}%
\pgfpathlineto{\pgfqpoint{3.318028in}{1.704656in}}%
\pgfpathlineto{\pgfqpoint{3.381463in}{1.697873in}}%
\pgfpathlineto{\pgfqpoint{3.440367in}{1.689287in}}%
\pgfpathlineto{\pgfqpoint{3.494740in}{1.679067in}}%
\pgfpathlineto{\pgfqpoint{3.544581in}{1.667365in}}%
\pgfpathlineto{\pgfqpoint{3.589892in}{1.654336in}}%
\pgfpathlineto{\pgfqpoint{3.630671in}{1.640169in}}%
\pgfpathlineto{\pgfqpoint{3.666920in}{1.625120in}}%
\pgfpathlineto{\pgfqpoint{3.698637in}{1.609572in}}%
\pgfpathlineto{\pgfqpoint{3.725824in}{1.594080in}}%
\pgfpathlineto{\pgfqpoint{3.753010in}{1.576268in}}%
\pgfpathlineto{\pgfqpoint{3.784728in}{1.552609in}}%
\pgfpathlineto{\pgfqpoint{3.820976in}{1.525343in}}%
\pgfpathlineto{\pgfqpoint{3.834569in}{1.517783in}}%
\pgfpathlineto{\pgfqpoint{3.843631in}{1.514580in}}%
\pgfpathlineto{\pgfqpoint{3.852693in}{1.513251in}}%
\pgfpathlineto{\pgfqpoint{3.861756in}{1.513951in}}%
\pgfpathlineto{\pgfqpoint{3.870818in}{1.516599in}}%
\pgfpathlineto{\pgfqpoint{3.884411in}{1.523556in}}%
\pgfpathlineto{\pgfqpoint{3.902535in}{1.536276in}}%
\pgfpathlineto{\pgfqpoint{3.965970in}{1.583696in}}%
\pgfpathlineto{\pgfqpoint{3.993156in}{1.600596in}}%
\pgfpathlineto{\pgfqpoint{4.024874in}{1.617541in}}%
\pgfpathlineto{\pgfqpoint{4.056591in}{1.631992in}}%
\pgfpathlineto{\pgfqpoint{4.092840in}{1.646042in}}%
\pgfpathlineto{\pgfqpoint{4.133619in}{1.659331in}}%
\pgfpathlineto{\pgfqpoint{4.178930in}{1.671597in}}%
\pgfpathlineto{\pgfqpoint{4.228771in}{1.682638in}}%
\pgfpathlineto{\pgfqpoint{4.283144in}{1.692285in}}%
\pgfpathlineto{\pgfqpoint{4.342048in}{1.700371in}}%
\pgfpathlineto{\pgfqpoint{4.405483in}{1.706710in}}%
\pgfpathlineto{\pgfqpoint{4.468917in}{1.710866in}}%
\pgfpathlineto{\pgfqpoint{4.536883in}{1.713082in}}%
\pgfpathlineto{\pgfqpoint{4.604849in}{1.713081in}}%
\pgfpathlineto{\pgfqpoint{4.672815in}{1.710867in}}%
\pgfpathlineto{\pgfqpoint{4.740781in}{1.706350in}}%
\pgfpathlineto{\pgfqpoint{4.804216in}{1.699888in}}%
\pgfpathlineto{\pgfqpoint{4.863120in}{1.691714in}}%
\pgfpathlineto{\pgfqpoint{4.917492in}{1.682031in}}%
\pgfpathlineto{\pgfqpoint{4.967334in}{1.671038in}}%
\pgfpathlineto{\pgfqpoint{5.012645in}{1.658957in}}%
\pgfpathlineto{\pgfqpoint{5.057955in}{1.644492in}}%
\pgfpathlineto{\pgfqpoint{5.098735in}{1.629061in}}%
\pgfpathlineto{\pgfqpoint{5.139514in}{1.611047in}}%
\pgfpathlineto{\pgfqpoint{5.184825in}{1.588256in}}%
\pgfpathlineto{\pgfqpoint{5.230135in}{1.565623in}}%
\pgfpathlineto{\pgfqpoint{5.248260in}{1.558659in}}%
\pgfpathlineto{\pgfqpoint{5.266384in}{1.554326in}}%
\pgfpathlineto{\pgfqpoint{5.279977in}{1.553313in}}%
\pgfpathlineto{\pgfqpoint{5.293570in}{1.554375in}}%
\pgfpathlineto{\pgfqpoint{5.307163in}{1.557384in}}%
\pgfpathlineto{\pgfqpoint{5.325288in}{1.563837in}}%
\pgfpathlineto{\pgfqpoint{5.352474in}{1.576690in}}%
\pgfpathlineto{\pgfqpoint{5.434033in}{1.617664in}}%
\pgfpathlineto{\pgfqpoint{5.474813in}{1.634854in}}%
\pgfpathlineto{\pgfqpoint{5.515592in}{1.649540in}}%
\pgfpathlineto{\pgfqpoint{5.529185in}{1.653934in}}%
\pgfpathlineto{\pgfqpoint{5.529185in}{1.653934in}}%
\pgfusepath{stroke}%
\end{pgfscope}%
\begin{pgfscope}%
\pgfpathrectangle{\pgfqpoint{0.776331in}{0.668778in}}{\pgfqpoint{4.979181in}{1.859128in}}%
\pgfusepath{clip}%
\pgfsetrectcap%
\pgfsetroundjoin%
\pgfsetlinewidth{1.505625pt}%
\definecolor{currentstroke}{rgb}{0.000000,0.694118,0.917647}%
\pgfsetstrokecolor{currentstroke}%
\pgfsetdash{}{0pt}%
\pgfpathmoveto{\pgfqpoint{1.028517in}{0.658778in}}%
\pgfpathlineto{\pgfqpoint{1.047968in}{2.044053in}}%
\pgfpathlineto{\pgfqpoint{1.093278in}{2.110140in}}%
\pgfpathlineto{\pgfqpoint{1.138589in}{2.147511in}}%
\pgfpathlineto{\pgfqpoint{1.183899in}{2.173151in}}%
\pgfpathlineto{\pgfqpoint{1.229210in}{2.192225in}}%
\pgfpathlineto{\pgfqpoint{1.274521in}{2.206962in}}%
\pgfpathlineto{\pgfqpoint{1.319831in}{2.218580in}}%
\pgfpathlineto{\pgfqpoint{1.365142in}{2.227891in}}%
\pgfpathlineto{\pgfqpoint{1.410452in}{2.235525in}}%
\pgfpathlineto{\pgfqpoint{1.455763in}{2.241997in}}%
\pgfpathlineto{\pgfqpoint{1.501074in}{2.247714in}}%
\pgfpathlineto{\pgfqpoint{1.546384in}{2.252962in}}%
\pgfpathlineto{\pgfqpoint{1.591695in}{2.257898in}}%
\pgfpathlineto{\pgfqpoint{1.637005in}{2.262563in}}%
\pgfpathlineto{\pgfqpoint{1.682316in}{2.266933in}}%
\pgfpathlineto{\pgfqpoint{1.727627in}{2.270970in}}%
\pgfpathlineto{\pgfqpoint{1.772937in}{2.274680in}}%
\pgfpathlineto{\pgfqpoint{1.818248in}{2.278130in}}%
\pgfpathlineto{\pgfqpoint{1.863558in}{2.281430in}}%
\pgfpathlineto{\pgfqpoint{1.908869in}{2.284686in}}%
\pgfpathlineto{\pgfqpoint{1.954180in}{2.287955in}}%
\pgfpathlineto{\pgfqpoint{1.999490in}{2.291234in}}%
\pgfpathlineto{\pgfqpoint{2.044801in}{2.294472in}}%
\pgfpathlineto{\pgfqpoint{2.090111in}{2.297608in}}%
\pgfpathlineto{\pgfqpoint{2.135422in}{2.300595in}}%
\pgfpathlineto{\pgfqpoint{2.180732in}{2.303422in}}%
\pgfpathlineto{\pgfqpoint{2.226043in}{2.306107in}}%
\pgfpathlineto{\pgfqpoint{2.271354in}{2.308684in}}%
\pgfpathlineto{\pgfqpoint{2.316664in}{2.311181in}}%
\pgfpathlineto{\pgfqpoint{2.361975in}{2.313602in}}%
\pgfpathlineto{\pgfqpoint{2.407285in}{2.315925in}}%
\pgfpathlineto{\pgfqpoint{2.452596in}{2.318115in}}%
\pgfpathlineto{\pgfqpoint{2.497907in}{2.320141in}}%
\pgfpathlineto{\pgfqpoint{2.543217in}{2.322000in}}%
\pgfpathlineto{\pgfqpoint{2.588528in}{2.323722in}}%
\pgfpathlineto{\pgfqpoint{2.633838in}{2.325355in}}%
\pgfpathlineto{\pgfqpoint{2.679149in}{2.326948in}}%
\pgfpathlineto{\pgfqpoint{2.724460in}{2.328531in}}%
\pgfpathlineto{\pgfqpoint{2.769770in}{2.330113in}}%
\pgfpathlineto{\pgfqpoint{2.815081in}{2.331687in}}%
\pgfpathlineto{\pgfqpoint{2.860391in}{2.333245in}}%
\pgfpathlineto{\pgfqpoint{2.905702in}{2.334781in}}%
\pgfpathlineto{\pgfqpoint{2.951013in}{2.336296in}}%
\pgfpathlineto{\pgfqpoint{2.996323in}{2.337792in}}%
\pgfpathlineto{\pgfqpoint{3.041634in}{2.339259in}}%
\pgfpathlineto{\pgfqpoint{3.086944in}{2.340677in}}%
\pgfpathlineto{\pgfqpoint{3.132255in}{2.342021in}}%
\pgfpathlineto{\pgfqpoint{3.177566in}{2.343268in}}%
\pgfpathlineto{\pgfqpoint{3.222876in}{2.344417in}}%
\pgfpathlineto{\pgfqpoint{3.268187in}{2.345487in}}%
\pgfpathlineto{\pgfqpoint{3.313497in}{2.346512in}}%
\pgfpathlineto{\pgfqpoint{3.358808in}{2.347525in}}%
\pgfpathlineto{\pgfqpoint{3.404119in}{2.348547in}}%
\pgfpathlineto{\pgfqpoint{3.449429in}{2.349586in}}%
\pgfpathlineto{\pgfqpoint{3.494740in}{2.350637in}}%
\pgfpathlineto{\pgfqpoint{3.540050in}{2.351692in}}%
\pgfpathlineto{\pgfqpoint{3.585361in}{2.352749in}}%
\pgfpathlineto{\pgfqpoint{3.630671in}{2.353808in}}%
\pgfpathlineto{\pgfqpoint{3.675982in}{2.354866in}}%
\pgfpathlineto{\pgfqpoint{3.721293in}{2.355917in}}%
\pgfpathlineto{\pgfqpoint{3.766603in}{2.356944in}}%
\pgfpathlineto{\pgfqpoint{3.811914in}{2.357925in}}%
\pgfpathlineto{\pgfqpoint{3.857224in}{2.358846in}}%
\pgfpathlineto{\pgfqpoint{3.902535in}{2.359703in}}%
\pgfpathlineto{\pgfqpoint{3.947846in}{2.360511in}}%
\pgfpathlineto{\pgfqpoint{3.993156in}{2.361293in}}%
\pgfpathlineto{\pgfqpoint{4.038467in}{2.362072in}}%
\pgfpathlineto{\pgfqpoint{4.083777in}{2.362862in}}%
\pgfpathlineto{\pgfqpoint{4.129088in}{2.363663in}}%
\pgfpathlineto{\pgfqpoint{4.174399in}{2.364466in}}%
\pgfpathlineto{\pgfqpoint{4.219709in}{2.365262in}}%
\pgfpathlineto{\pgfqpoint{4.265020in}{2.366044in}}%
\pgfpathlineto{\pgfqpoint{4.310330in}{2.366816in}}%
\pgfpathlineto{\pgfqpoint{4.355641in}{2.367584in}}%
\pgfpathlineto{\pgfqpoint{4.400952in}{2.368354in}}%
\pgfpathlineto{\pgfqpoint{4.446262in}{2.369128in}}%
\pgfpathlineto{\pgfqpoint{4.491573in}{2.369898in}}%
\pgfpathlineto{\pgfqpoint{4.536883in}{2.370646in}}%
\pgfpathlineto{\pgfqpoint{4.582194in}{2.371357in}}%
\pgfpathlineto{\pgfqpoint{4.627505in}{2.372020in}}%
\pgfpathlineto{\pgfqpoint{4.672815in}{2.372637in}}%
\pgfpathlineto{\pgfqpoint{4.718126in}{2.373225in}}%
\pgfpathlineto{\pgfqpoint{4.763436in}{2.373803in}}%
\pgfpathlineto{\pgfqpoint{4.808747in}{2.374388in}}%
\pgfpathlineto{\pgfqpoint{4.854058in}{2.374987in}}%
\pgfpathlineto{\pgfqpoint{4.899368in}{2.375598in}}%
\pgfpathlineto{\pgfqpoint{4.944679in}{2.376213in}}%
\pgfpathlineto{\pgfqpoint{4.989989in}{2.376825in}}%
\pgfpathlineto{\pgfqpoint{5.035300in}{2.377431in}}%
\pgfpathlineto{\pgfqpoint{5.080611in}{2.378034in}}%
\pgfpathlineto{\pgfqpoint{5.125921in}{2.378640in}}%
\pgfpathlineto{\pgfqpoint{5.171232in}{2.379251in}}%
\pgfpathlineto{\pgfqpoint{5.216542in}{2.379867in}}%
\pgfpathlineto{\pgfqpoint{5.261853in}{2.380477in}}%
\pgfpathlineto{\pgfqpoint{5.307163in}{2.381064in}}%
\pgfpathlineto{\pgfqpoint{5.352474in}{2.381618in}}%
\pgfpathlineto{\pgfqpoint{5.397785in}{2.382131in}}%
\pgfpathlineto{\pgfqpoint{5.443095in}{2.382611in}}%
\pgfpathlineto{\pgfqpoint{5.488406in}{2.383071in}}%
\pgfusepath{stroke}%
\end{pgfscope}%
\begin{pgfscope}%
\pgfpathrectangle{\pgfqpoint{0.776331in}{0.668778in}}{\pgfqpoint{4.979181in}{1.859128in}}%
\pgfusepath{clip}%
\pgfsetrectcap%
\pgfsetroundjoin%
\pgfsetlinewidth{1.505625pt}%
\definecolor{currentstroke}{rgb}{0.937255,0.486275,0.000000}%
\pgfsetstrokecolor{currentstroke}%
\pgfsetdash{}{0pt}%
\pgfpathmoveto{\pgfqpoint{1.030622in}{0.658778in}}%
\pgfpathlineto{\pgfqpoint{1.047968in}{1.801151in}}%
\pgfpathlineto{\pgfqpoint{1.093278in}{1.870892in}}%
\pgfpathlineto{\pgfqpoint{1.138589in}{1.911553in}}%
\pgfpathlineto{\pgfqpoint{1.183899in}{1.940284in}}%
\pgfpathlineto{\pgfqpoint{1.229210in}{1.962458in}}%
\pgfpathlineto{\pgfqpoint{1.274521in}{1.980467in}}%
\pgfpathlineto{\pgfqpoint{1.319831in}{1.995586in}}%
\pgfpathlineto{\pgfqpoint{1.365142in}{2.008575in}}%
\pgfpathlineto{\pgfqpoint{1.410452in}{2.019923in}}%
\pgfpathlineto{\pgfqpoint{1.455763in}{2.029962in}}%
\pgfpathlineto{\pgfqpoint{1.501074in}{2.038930in}}%
\pgfpathlineto{\pgfqpoint{1.546384in}{2.047001in}}%
\pgfpathlineto{\pgfqpoint{1.591695in}{2.054307in}}%
\pgfpathlineto{\pgfqpoint{1.637005in}{2.060953in}}%
\pgfpathlineto{\pgfqpoint{1.682316in}{2.067019in}}%
\pgfpathlineto{\pgfqpoint{1.727627in}{2.072584in}}%
\pgfpathlineto{\pgfqpoint{1.772937in}{2.077834in}}%
\pgfpathlineto{\pgfqpoint{1.818248in}{2.082823in}}%
\pgfpathlineto{\pgfqpoint{1.863558in}{2.087553in}}%
\pgfpathlineto{\pgfqpoint{1.908869in}{2.092031in}}%
\pgfpathlineto{\pgfqpoint{1.954180in}{2.096268in}}%
\pgfpathlineto{\pgfqpoint{1.999490in}{2.100275in}}%
\pgfpathlineto{\pgfqpoint{2.044801in}{2.104064in}}%
\pgfpathlineto{\pgfqpoint{2.090111in}{2.107645in}}%
\pgfpathlineto{\pgfqpoint{2.135422in}{2.111029in}}%
\pgfpathlineto{\pgfqpoint{2.180732in}{2.114224in}}%
\pgfpathlineto{\pgfqpoint{2.226043in}{2.117239in}}%
\pgfpathlineto{\pgfqpoint{2.271354in}{2.120081in}}%
\pgfpathlineto{\pgfqpoint{2.316664in}{2.122755in}}%
\pgfpathlineto{\pgfqpoint{2.361975in}{2.125268in}}%
\pgfpathlineto{\pgfqpoint{2.407285in}{2.127626in}}%
\pgfpathlineto{\pgfqpoint{2.452596in}{2.129855in}}%
\pgfpathlineto{\pgfqpoint{2.497907in}{2.132048in}}%
\pgfpathlineto{\pgfqpoint{2.543217in}{2.134208in}}%
\pgfpathlineto{\pgfqpoint{2.588528in}{2.136329in}}%
\pgfpathlineto{\pgfqpoint{2.633838in}{2.138407in}}%
\pgfpathlineto{\pgfqpoint{2.679149in}{2.140439in}}%
\pgfpathlineto{\pgfqpoint{2.724460in}{2.142423in}}%
\pgfpathlineto{\pgfqpoint{2.769770in}{2.144359in}}%
\pgfpathlineto{\pgfqpoint{2.815081in}{2.146247in}}%
\pgfpathlineto{\pgfqpoint{2.860391in}{2.148085in}}%
\pgfpathlineto{\pgfqpoint{2.905702in}{2.149875in}}%
\pgfpathlineto{\pgfqpoint{2.951013in}{2.151615in}}%
\pgfpathlineto{\pgfqpoint{2.996323in}{2.153304in}}%
\pgfpathlineto{\pgfqpoint{3.041634in}{2.154943in}}%
\pgfpathlineto{\pgfqpoint{3.086944in}{2.156530in}}%
\pgfpathlineto{\pgfqpoint{3.132255in}{2.158063in}}%
\pgfpathlineto{\pgfqpoint{3.177566in}{2.159585in}}%
\pgfpathlineto{\pgfqpoint{3.222876in}{2.161159in}}%
\pgfpathlineto{\pgfqpoint{3.268187in}{2.162773in}}%
\pgfpathlineto{\pgfqpoint{3.313497in}{2.164418in}}%
\pgfpathlineto{\pgfqpoint{3.358808in}{2.166084in}}%
\pgfpathlineto{\pgfqpoint{3.404119in}{2.167762in}}%
\pgfpathlineto{\pgfqpoint{3.449429in}{2.169443in}}%
\pgfpathlineto{\pgfqpoint{3.494740in}{2.171120in}}%
\pgfpathlineto{\pgfqpoint{3.540050in}{2.172785in}}%
\pgfpathlineto{\pgfqpoint{3.585361in}{2.174433in}}%
\pgfpathlineto{\pgfqpoint{3.630671in}{2.176057in}}%
\pgfpathlineto{\pgfqpoint{3.675982in}{2.177653in}}%
\pgfpathlineto{\pgfqpoint{3.721293in}{2.179215in}}%
\pgfpathlineto{\pgfqpoint{3.766603in}{2.180741in}}%
\pgfpathlineto{\pgfqpoint{3.811914in}{2.182226in}}%
\pgfpathlineto{\pgfqpoint{3.857224in}{2.183669in}}%
\pgfpathlineto{\pgfqpoint{3.902535in}{2.185123in}}%
\pgfpathlineto{\pgfqpoint{3.947846in}{2.186610in}}%
\pgfpathlineto{\pgfqpoint{3.993156in}{2.188121in}}%
\pgfpathlineto{\pgfqpoint{4.038467in}{2.189648in}}%
\pgfpathlineto{\pgfqpoint{4.083777in}{2.191181in}}%
\pgfpathlineto{\pgfqpoint{4.129088in}{2.192714in}}%
\pgfpathlineto{\pgfqpoint{4.174399in}{2.194238in}}%
\pgfpathlineto{\pgfqpoint{4.219709in}{2.195747in}}%
\pgfpathlineto{\pgfqpoint{4.265020in}{2.197235in}}%
\pgfpathlineto{\pgfqpoint{4.310330in}{2.198696in}}%
\pgfpathlineto{\pgfqpoint{4.355641in}{2.200125in}}%
\pgfpathlineto{\pgfqpoint{4.400952in}{2.201518in}}%
\pgfpathlineto{\pgfqpoint{4.446262in}{2.202871in}}%
\pgfpathlineto{\pgfqpoint{4.491573in}{2.204179in}}%
\pgfpathlineto{\pgfqpoint{4.536883in}{2.205442in}}%
\pgfpathlineto{\pgfqpoint{4.582194in}{2.206663in}}%
\pgfpathlineto{\pgfqpoint{4.627505in}{2.207888in}}%
\pgfpathlineto{\pgfqpoint{4.672815in}{2.209121in}}%
\pgfpathlineto{\pgfqpoint{4.718126in}{2.210358in}}%
\pgfpathlineto{\pgfqpoint{4.763436in}{2.211592in}}%
\pgfpathlineto{\pgfqpoint{4.808747in}{2.212819in}}%
\pgfpathlineto{\pgfqpoint{4.854058in}{2.214034in}}%
\pgfpathlineto{\pgfqpoint{4.899368in}{2.215234in}}%
\pgfpathlineto{\pgfqpoint{4.944679in}{2.216414in}}%
\pgfpathlineto{\pgfqpoint{4.989989in}{2.217571in}}%
\pgfpathlineto{\pgfqpoint{5.035300in}{2.218701in}}%
\pgfpathlineto{\pgfqpoint{5.080611in}{2.219802in}}%
\pgfpathlineto{\pgfqpoint{5.125921in}{2.220870in}}%
\pgfpathlineto{\pgfqpoint{5.171232in}{2.221903in}}%
\pgfpathlineto{\pgfqpoint{5.216542in}{2.222899in}}%
\pgfpathlineto{\pgfqpoint{5.261853in}{2.223856in}}%
\pgfpathlineto{\pgfqpoint{5.307163in}{2.224787in}}%
\pgfpathlineto{\pgfqpoint{5.352474in}{2.225724in}}%
\pgfpathlineto{\pgfqpoint{5.397785in}{2.226667in}}%
\pgfpathlineto{\pgfqpoint{5.443095in}{2.227611in}}%
\pgfpathlineto{\pgfqpoint{5.488406in}{2.228553in}}%
\pgfusepath{stroke}%
\end{pgfscope}%
\begin{pgfscope}%
\pgfsetrectcap%
\pgfsetmiterjoin%
\pgfsetlinewidth{0.803000pt}%
\definecolor{currentstroke}{rgb}{0.000000,0.000000,0.000000}%
\pgfsetstrokecolor{currentstroke}%
\pgfsetdash{}{0pt}%
\pgfpathmoveto{\pgfqpoint{0.776331in}{0.668778in}}%
\pgfpathlineto{\pgfqpoint{0.776331in}{2.527906in}}%
\pgfusepath{stroke}%
\end{pgfscope}%
\begin{pgfscope}%
\pgfsetrectcap%
\pgfsetmiterjoin%
\pgfsetlinewidth{0.803000pt}%
\definecolor{currentstroke}{rgb}{0.000000,0.000000,0.000000}%
\pgfsetstrokecolor{currentstroke}%
\pgfsetdash{}{0pt}%
\pgfpathmoveto{\pgfqpoint{5.755512in}{0.668778in}}%
\pgfpathlineto{\pgfqpoint{5.755512in}{2.527906in}}%
\pgfusepath{stroke}%
\end{pgfscope}%
\begin{pgfscope}%
\pgfsetrectcap%
\pgfsetmiterjoin%
\pgfsetlinewidth{0.803000pt}%
\definecolor{currentstroke}{rgb}{0.000000,0.000000,0.000000}%
\pgfsetstrokecolor{currentstroke}%
\pgfsetdash{}{0pt}%
\pgfpathmoveto{\pgfqpoint{0.776331in}{0.668778in}}%
\pgfpathlineto{\pgfqpoint{5.755512in}{0.668778in}}%
\pgfusepath{stroke}%
\end{pgfscope}%
\begin{pgfscope}%
\pgfsetrectcap%
\pgfsetmiterjoin%
\pgfsetlinewidth{0.803000pt}%
\definecolor{currentstroke}{rgb}{0.000000,0.000000,0.000000}%
\pgfsetstrokecolor{currentstroke}%
\pgfsetdash{}{0pt}%
\pgfpathmoveto{\pgfqpoint{0.776331in}{2.527906in}}%
\pgfpathlineto{\pgfqpoint{5.755512in}{2.527906in}}%
\pgfusepath{stroke}%
\end{pgfscope}%
\begin{pgfscope}%
\pgfsetbuttcap%
\pgfsetmiterjoin%
\definecolor{currentfill}{rgb}{1.000000,1.000000,1.000000}%
\pgfsetfillcolor{currentfill}%
\pgfsetfillopacity{0.800000}%
\pgfsetlinewidth{1.003750pt}%
\definecolor{currentstroke}{rgb}{0.800000,0.800000,0.800000}%
\pgfsetstrokecolor{currentstroke}%
\pgfsetstrokeopacity{0.800000}%
\pgfsetdash{}{0pt}%
\pgfpathmoveto{\pgfqpoint{3.816140in}{0.738222in}}%
\pgfpathlineto{\pgfqpoint{5.658290in}{0.738222in}}%
\pgfpathquadraticcurveto{\pgfqpoint{5.686067in}{0.738222in}}{\pgfqpoint{5.686067in}{0.766000in}}%
\pgfpathlineto{\pgfqpoint{5.686067in}{1.363683in}}%
\pgfpathquadraticcurveto{\pgfqpoint{5.686067in}{1.391461in}}{\pgfqpoint{5.658290in}{1.391461in}}%
\pgfpathlineto{\pgfqpoint{3.816140in}{1.391461in}}%
\pgfpathquadraticcurveto{\pgfqpoint{3.788362in}{1.391461in}}{\pgfqpoint{3.788362in}{1.363683in}}%
\pgfpathlineto{\pgfqpoint{3.788362in}{0.766000in}}%
\pgfpathquadraticcurveto{\pgfqpoint{3.788362in}{0.738222in}}{\pgfqpoint{3.816140in}{0.738222in}}%
\pgfpathlineto{\pgfqpoint{3.816140in}{0.738222in}}%
\pgfpathclose%
\pgfusepath{stroke,fill}%
\end{pgfscope}%
\begin{pgfscope}%
\pgfsetrectcap%
\pgfsetroundjoin%
\pgfsetlinewidth{1.505625pt}%
\definecolor{currentstroke}{rgb}{0.776471,0.827451,0.145098}%
\pgfsetstrokecolor{currentstroke}%
\pgfsetdash{}{0pt}%
\pgfpathmoveto{\pgfqpoint{3.843918in}{1.278993in}}%
\pgfpathlineto{\pgfqpoint{3.982807in}{1.278993in}}%
\pgfpathlineto{\pgfqpoint{4.121696in}{1.278993in}}%
\pgfusepath{stroke}%
\end{pgfscope}%
\begin{pgfscope}%
\definecolor{textcolor}{rgb}{0.000000,0.000000,0.000000}%
\pgfsetstrokecolor{textcolor}%
\pgfsetfillcolor{textcolor}%
\pgftext[x=4.232807in,y=1.230382in,left,base]{\color{textcolor}{\sffamily\fontsize{10.000000}{12.000000}\selectfont\catcode`\^=\active\def^{\ifmmode\sp\else\^{}\fi}\catcode`\%=\active\def
\end{pgfscope}%
\begin{pgfscope}%
\pgfsetrectcap%
\pgfsetroundjoin%
\pgfsetlinewidth{1.505625pt}%
\definecolor{currentstroke}{rgb}{0.000000,0.694118,0.917647}%
\pgfsetstrokecolor{currentstroke}%
\pgfsetdash{}{0pt}%
\pgfpathmoveto{\pgfqpoint{3.843918in}{1.075136in}}%
\pgfpathlineto{\pgfqpoint{3.982807in}{1.075136in}}%
\pgfpathlineto{\pgfqpoint{4.121696in}{1.075136in}}%
\pgfusepath{stroke}%
\end{pgfscope}%
\begin{pgfscope}%
\definecolor{textcolor}{rgb}{0.000000,0.000000,0.000000}%
\pgfsetstrokecolor{textcolor}%
\pgfsetfillcolor{textcolor}%
\pgftext[x=4.232807in,y=1.026525in,left,base]{\color{textcolor}{\sffamily\fontsize{10.000000}{12.000000}\selectfont\catcode`\^=\active\def^{\ifmmode\sp\else\^{}\fi}\catcode`\%=\active\def
\end{pgfscope}%
\begin{pgfscope}%
\pgfsetrectcap%
\pgfsetroundjoin%
\pgfsetlinewidth{1.505625pt}%
\definecolor{currentstroke}{rgb}{0.937255,0.486275,0.000000}%
\pgfsetstrokecolor{currentstroke}%
\pgfsetdash{}{0pt}%
\pgfpathmoveto{\pgfqpoint{3.843918in}{0.871279in}}%
\pgfpathlineto{\pgfqpoint{3.982807in}{0.871279in}}%
\pgfpathlineto{\pgfqpoint{4.121696in}{0.871279in}}%
\pgfusepath{stroke}%
\end{pgfscope}%
\begin{pgfscope}%
\definecolor{textcolor}{rgb}{0.000000,0.000000,0.000000}%
\pgfsetstrokecolor{textcolor}%
\pgfsetfillcolor{textcolor}%
\pgftext[x=4.232807in,y=0.822667in,left,base]{\color{textcolor}{\sffamily\fontsize{10.000000}{12.000000}\selectfont\catcode`\^=\active\def^{\ifmmode\sp\else\^{}\fi}\catcode`\%=\active\def
\end{pgfscope}%
\end{pgfpicture}%
\makeatother%
\endgroup%

%% file: Graphics/RotatingGrid_convergenceRates.pgf
\begingroup%
\makeatletter%
\begin{pgfpicture}%
\pgfpathrectangle{\pgfpointorigin}{\pgfqpoint{7.086614in}{3.937008in}}%
\pgfusepath{use as bounding box, clip}%
\begin{pgfscope}%
\pgfsetbuttcap%
\pgfsetmiterjoin%
\definecolor{currentfill}{rgb}{1.000000,1.000000,1.000000}%
\pgfsetfillcolor{currentfill}%
\pgfsetlinewidth{0.000000pt}%
\definecolor{currentstroke}{rgb}{1.000000,1.000000,1.000000}%
\pgfsetstrokecolor{currentstroke}%
\pgfsetdash{}{0pt}%
\pgfpathmoveto{\pgfqpoint{0.000000in}{0.000000in}}%
\pgfpathlineto{\pgfqpoint{7.086614in}{0.000000in}}%
\pgfpathlineto{\pgfqpoint{7.086614in}{3.937008in}}%
\pgfpathlineto{\pgfqpoint{0.000000in}{3.937008in}}%
\pgfpathlineto{\pgfqpoint{0.000000in}{0.000000in}}%
\pgfpathclose%
\pgfusepath{fill}%
\end{pgfscope}%
\begin{pgfscope}%
\pgfsetbuttcap%
\pgfsetmiterjoin%
\definecolor{currentfill}{rgb}{1.000000,1.000000,1.000000}%
\pgfsetfillcolor{currentfill}%
\pgfsetlinewidth{0.000000pt}%
\definecolor{currentstroke}{rgb}{0.000000,0.000000,0.000000}%
\pgfsetstrokecolor{currentstroke}%
\pgfsetstrokeopacity{0.000000}%
\pgfsetdash{}{0pt}%
\pgfpathmoveto{\pgfqpoint{0.885827in}{0.433071in}}%
\pgfpathlineto{\pgfqpoint{6.377953in}{0.433071in}}%
\pgfpathlineto{\pgfqpoint{6.377953in}{3.464567in}}%
\pgfpathlineto{\pgfqpoint{0.885827in}{3.464567in}}%
\pgfpathlineto{\pgfqpoint{0.885827in}{0.433071in}}%
\pgfpathclose%
\pgfusepath{fill}%
\end{pgfscope}%
\begin{pgfscope}%
\pgfsetbuttcap%
\pgfsetroundjoin%
\definecolor{currentfill}{rgb}{0.000000,0.000000,0.000000}%
\pgfsetfillcolor{currentfill}%
\pgfsetlinewidth{0.803000pt}%
\definecolor{currentstroke}{rgb}{0.000000,0.000000,0.000000}%
\pgfsetstrokecolor{currentstroke}%
\pgfsetdash{}{0pt}%
\pgfsys@defobject{currentmarker}{\pgfqpoint{0.000000in}{-0.048611in}}{\pgfqpoint{0.000000in}{0.000000in}}{%
\pgfpathmoveto{\pgfqpoint{0.000000in}{0.000000in}}%
\pgfpathlineto{\pgfqpoint{0.000000in}{-0.048611in}}%
\pgfusepath{stroke,fill}%
}%
\begin{pgfscope}%
\pgfsys@transformshift{0.885827in}{0.433071in}%
\pgfsys@useobject{currentmarker}{}%
\end{pgfscope}%
\end{pgfscope}%
\begin{pgfscope}%
\definecolor{textcolor}{rgb}{0.000000,0.000000,0.000000}%
\pgfsetstrokecolor{textcolor}%
\pgfsetfillcolor{textcolor}%
\pgftext[x=0.885827in,y=0.335849in,,top]{\color{textcolor}{\sffamily\fontsize{10.000000}{12.000000}\selectfont\catcode`\^=\active\def^{\ifmmode\sp\else\^{}\fi}\catcode`\%=\active\def
\end{pgfscope}%
\begin{pgfscope}%
\pgfsetbuttcap%
\pgfsetroundjoin%
\definecolor{currentfill}{rgb}{0.000000,0.000000,0.000000}%
\pgfsetfillcolor{currentfill}%
\pgfsetlinewidth{0.803000pt}%
\definecolor{currentstroke}{rgb}{0.000000,0.000000,0.000000}%
\pgfsetstrokecolor{currentstroke}%
\pgfsetdash{}{0pt}%
\pgfsys@defobject{currentmarker}{\pgfqpoint{0.000000in}{-0.048611in}}{\pgfqpoint{0.000000in}{0.000000in}}{%
\pgfpathmoveto{\pgfqpoint{0.000000in}{0.000000in}}%
\pgfpathlineto{\pgfqpoint{0.000000in}{-0.048611in}}%
\pgfusepath{stroke,fill}%
}%
\begin{pgfscope}%
\pgfsys@transformshift{4.814555in}{0.433071in}%
\pgfsys@useobject{currentmarker}{}%
\end{pgfscope}%
\end{pgfscope}%
\begin{pgfscope}%
\definecolor{textcolor}{rgb}{0.000000,0.000000,0.000000}%
\pgfsetstrokecolor{textcolor}%
\pgfsetfillcolor{textcolor}%
\pgftext[x=4.814555in,y=0.335849in,,top]{\color{textcolor}{\sffamily\fontsize{10.000000}{12.000000}\selectfont\catcode`\^=\active\def^{\ifmmode\sp\else\^{}\fi}\catcode`\%=\active\def
\end{pgfscope}%
\begin{pgfscope}%
\pgfsetbuttcap%
\pgfsetroundjoin%
\definecolor{currentfill}{rgb}{0.000000,0.000000,0.000000}%
\pgfsetfillcolor{currentfill}%
\pgfsetlinewidth{0.602250pt}%
\definecolor{currentstroke}{rgb}{0.000000,0.000000,0.000000}%
\pgfsetstrokecolor{currentstroke}%
\pgfsetdash{}{0pt}%
\pgfsys@defobject{currentmarker}{\pgfqpoint{0.000000in}{-0.027778in}}{\pgfqpoint{0.000000in}{0.000000in}}{%
\pgfpathmoveto{\pgfqpoint{0.000000in}{0.000000in}}%
\pgfpathlineto{\pgfqpoint{0.000000in}{-0.027778in}}%
\pgfusepath{stroke,fill}%
}%
\begin{pgfscope}%
\pgfsys@transformshift{2.068492in}{0.433071in}%
\pgfsys@useobject{currentmarker}{}%
\end{pgfscope}%
\end{pgfscope}%
\begin{pgfscope}%
\pgfsetbuttcap%
\pgfsetroundjoin%
\definecolor{currentfill}{rgb}{0.000000,0.000000,0.000000}%
\pgfsetfillcolor{currentfill}%
\pgfsetlinewidth{0.602250pt}%
\definecolor{currentstroke}{rgb}{0.000000,0.000000,0.000000}%
\pgfsetstrokecolor{currentstroke}%
\pgfsetdash{}{0pt}%
\pgfsys@defobject{currentmarker}{\pgfqpoint{0.000000in}{-0.027778in}}{\pgfqpoint{0.000000in}{0.000000in}}{%
\pgfpathmoveto{\pgfqpoint{0.000000in}{0.000000in}}%
\pgfpathlineto{\pgfqpoint{0.000000in}{-0.027778in}}%
\pgfusepath{stroke,fill}%
}%
\begin{pgfscope}%
\pgfsys@transformshift{2.760306in}{0.433071in}%
\pgfsys@useobject{currentmarker}{}%
\end{pgfscope}%
\end{pgfscope}%
\begin{pgfscope}%
\pgfsetbuttcap%
\pgfsetroundjoin%
\definecolor{currentfill}{rgb}{0.000000,0.000000,0.000000}%
\pgfsetfillcolor{currentfill}%
\pgfsetlinewidth{0.602250pt}%
\definecolor{currentstroke}{rgb}{0.000000,0.000000,0.000000}%
\pgfsetstrokecolor{currentstroke}%
\pgfsetdash{}{0pt}%
\pgfsys@defobject{currentmarker}{\pgfqpoint{0.000000in}{-0.027778in}}{\pgfqpoint{0.000000in}{0.000000in}}{%
\pgfpathmoveto{\pgfqpoint{0.000000in}{0.000000in}}%
\pgfpathlineto{\pgfqpoint{0.000000in}{-0.027778in}}%
\pgfusepath{stroke,fill}%
}%
\begin{pgfscope}%
\pgfsys@transformshift{3.251157in}{0.433071in}%
\pgfsys@useobject{currentmarker}{}%
\end{pgfscope}%
\end{pgfscope}%
\begin{pgfscope}%
\pgfsetbuttcap%
\pgfsetroundjoin%
\definecolor{currentfill}{rgb}{0.000000,0.000000,0.000000}%
\pgfsetfillcolor{currentfill}%
\pgfsetlinewidth{0.602250pt}%
\definecolor{currentstroke}{rgb}{0.000000,0.000000,0.000000}%
\pgfsetstrokecolor{currentstroke}%
\pgfsetdash{}{0pt}%
\pgfsys@defobject{currentmarker}{\pgfqpoint{0.000000in}{-0.027778in}}{\pgfqpoint{0.000000in}{0.000000in}}{%
\pgfpathmoveto{\pgfqpoint{0.000000in}{0.000000in}}%
\pgfpathlineto{\pgfqpoint{0.000000in}{-0.027778in}}%
\pgfusepath{stroke,fill}%
}%
\begin{pgfscope}%
\pgfsys@transformshift{3.631890in}{0.433071in}%
\pgfsys@useobject{currentmarker}{}%
\end{pgfscope}%
\end{pgfscope}%
\begin{pgfscope}%
\pgfsetbuttcap%
\pgfsetroundjoin%
\definecolor{currentfill}{rgb}{0.000000,0.000000,0.000000}%
\pgfsetfillcolor{currentfill}%
\pgfsetlinewidth{0.602250pt}%
\definecolor{currentstroke}{rgb}{0.000000,0.000000,0.000000}%
\pgfsetstrokecolor{currentstroke}%
\pgfsetdash{}{0pt}%
\pgfsys@defobject{currentmarker}{\pgfqpoint{0.000000in}{-0.027778in}}{\pgfqpoint{0.000000in}{0.000000in}}{%
\pgfpathmoveto{\pgfqpoint{0.000000in}{0.000000in}}%
\pgfpathlineto{\pgfqpoint{0.000000in}{-0.027778in}}%
\pgfusepath{stroke,fill}%
}%
\begin{pgfscope}%
\pgfsys@transformshift{3.942971in}{0.433071in}%
\pgfsys@useobject{currentmarker}{}%
\end{pgfscope}%
\end{pgfscope}%
\begin{pgfscope}%
\pgfsetbuttcap%
\pgfsetroundjoin%
\definecolor{currentfill}{rgb}{0.000000,0.000000,0.000000}%
\pgfsetfillcolor{currentfill}%
\pgfsetlinewidth{0.602250pt}%
\definecolor{currentstroke}{rgb}{0.000000,0.000000,0.000000}%
\pgfsetstrokecolor{currentstroke}%
\pgfsetdash{}{0pt}%
\pgfsys@defobject{currentmarker}{\pgfqpoint{0.000000in}{-0.027778in}}{\pgfqpoint{0.000000in}{0.000000in}}{%
\pgfpathmoveto{\pgfqpoint{0.000000in}{0.000000in}}%
\pgfpathlineto{\pgfqpoint{0.000000in}{-0.027778in}}%
\pgfusepath{stroke,fill}%
}%
\begin{pgfscope}%
\pgfsys@transformshift{4.205987in}{0.433071in}%
\pgfsys@useobject{currentmarker}{}%
\end{pgfscope}%
\end{pgfscope}%
\begin{pgfscope}%
\pgfsetbuttcap%
\pgfsetroundjoin%
\definecolor{currentfill}{rgb}{0.000000,0.000000,0.000000}%
\pgfsetfillcolor{currentfill}%
\pgfsetlinewidth{0.602250pt}%
\definecolor{currentstroke}{rgb}{0.000000,0.000000,0.000000}%
\pgfsetstrokecolor{currentstroke}%
\pgfsetdash{}{0pt}%
\pgfsys@defobject{currentmarker}{\pgfqpoint{0.000000in}{-0.027778in}}{\pgfqpoint{0.000000in}{0.000000in}}{%
\pgfpathmoveto{\pgfqpoint{0.000000in}{0.000000in}}%
\pgfpathlineto{\pgfqpoint{0.000000in}{-0.027778in}}%
\pgfusepath{stroke,fill}%
}%
\begin{pgfscope}%
\pgfsys@transformshift{4.433822in}{0.433071in}%
\pgfsys@useobject{currentmarker}{}%
\end{pgfscope}%
\end{pgfscope}%
\begin{pgfscope}%
\pgfsetbuttcap%
\pgfsetroundjoin%
\definecolor{currentfill}{rgb}{0.000000,0.000000,0.000000}%
\pgfsetfillcolor{currentfill}%
\pgfsetlinewidth{0.602250pt}%
\definecolor{currentstroke}{rgb}{0.000000,0.000000,0.000000}%
\pgfsetstrokecolor{currentstroke}%
\pgfsetdash{}{0pt}%
\pgfsys@defobject{currentmarker}{\pgfqpoint{0.000000in}{-0.027778in}}{\pgfqpoint{0.000000in}{0.000000in}}{%
\pgfpathmoveto{\pgfqpoint{0.000000in}{0.000000in}}%
\pgfpathlineto{\pgfqpoint{0.000000in}{-0.027778in}}%
\pgfusepath{stroke,fill}%
}%
\begin{pgfscope}%
\pgfsys@transformshift{4.634786in}{0.433071in}%
\pgfsys@useobject{currentmarker}{}%
\end{pgfscope}%
\end{pgfscope}%
\begin{pgfscope}%
\pgfsetbuttcap%
\pgfsetroundjoin%
\definecolor{currentfill}{rgb}{0.000000,0.000000,0.000000}%
\pgfsetfillcolor{currentfill}%
\pgfsetlinewidth{0.602250pt}%
\definecolor{currentstroke}{rgb}{0.000000,0.000000,0.000000}%
\pgfsetstrokecolor{currentstroke}%
\pgfsetdash{}{0pt}%
\pgfsys@defobject{currentmarker}{\pgfqpoint{0.000000in}{-0.027778in}}{\pgfqpoint{0.000000in}{0.000000in}}{%
\pgfpathmoveto{\pgfqpoint{0.000000in}{0.000000in}}%
\pgfpathlineto{\pgfqpoint{0.000000in}{-0.027778in}}%
\pgfusepath{stroke,fill}%
}%
\begin{pgfscope}%
\pgfsys@transformshift{5.997220in}{0.433071in}%
\pgfsys@useobject{currentmarker}{}%
\end{pgfscope}%
\end{pgfscope}%
\begin{pgfscope}%
\definecolor{textcolor}{rgb}{0.000000,0.000000,0.000000}%
\pgfsetstrokecolor{textcolor}%
\pgfsetfillcolor{textcolor}%
\pgftext[x=3.631890in,y=0.145880in,,top]{\color{textcolor}{\sffamily\fontsize{10.000000}{12.000000}\selectfont\catcode`\^=\active\def^{\ifmmode\sp\else\^{}\fi}\catcode`\%=\active\def
\end{pgfscope}%
\begin{pgfscope}%
\pgfsetbuttcap%
\pgfsetroundjoin%
\definecolor{currentfill}{rgb}{0.000000,0.000000,0.000000}%
\pgfsetfillcolor{currentfill}%
\pgfsetlinewidth{0.803000pt}%
\definecolor{currentstroke}{rgb}{0.000000,0.000000,0.000000}%
\pgfsetstrokecolor{currentstroke}%
\pgfsetdash{}{0pt}%
\pgfsys@defobject{currentmarker}{\pgfqpoint{-0.048611in}{0.000000in}}{\pgfqpoint{-0.000000in}{0.000000in}}{%
\pgfpathmoveto{\pgfqpoint{-0.000000in}{0.000000in}}%
\pgfpathlineto{\pgfqpoint{-0.048611in}{0.000000in}}%
\pgfusepath{stroke,fill}%
}%
\begin{pgfscope}%
\pgfsys@transformshift{0.885827in}{0.785274in}%
\pgfsys@useobject{currentmarker}{}%
\end{pgfscope}%
\end{pgfscope}%
\begin{pgfscope}%
\definecolor{textcolor}{rgb}{0.000000,0.000000,0.000000}%
\pgfsetstrokecolor{textcolor}%
\pgfsetfillcolor{textcolor}%
\pgftext[x=0.445239in, y=0.732512in, left, base]{\color{textcolor}{\sffamily\fontsize{10.000000}{12.000000}\selectfont\catcode`\^=\active\def^{\ifmmode\sp\else\^{}\fi}\catcode`\%=\active\def
\end{pgfscope}%
\begin{pgfscope}%
\pgfsetbuttcap%
\pgfsetroundjoin%
\definecolor{currentfill}{rgb}{0.000000,0.000000,0.000000}%
\pgfsetfillcolor{currentfill}%
\pgfsetlinewidth{0.803000pt}%
\definecolor{currentstroke}{rgb}{0.000000,0.000000,0.000000}%
\pgfsetstrokecolor{currentstroke}%
\pgfsetdash{}{0pt}%
\pgfsys@defobject{currentmarker}{\pgfqpoint{-0.048611in}{0.000000in}}{\pgfqpoint{-0.000000in}{0.000000in}}{%
\pgfpathmoveto{\pgfqpoint{-0.000000in}{0.000000in}}%
\pgfpathlineto{\pgfqpoint{-0.048611in}{0.000000in}}%
\pgfusepath{stroke,fill}%
}%
\begin{pgfscope}%
\pgfsys@transformshift{0.885827in}{1.353161in}%
\pgfsys@useobject{currentmarker}{}%
\end{pgfscope}%
\end{pgfscope}%
\begin{pgfscope}%
\definecolor{textcolor}{rgb}{0.000000,0.000000,0.000000}%
\pgfsetstrokecolor{textcolor}%
\pgfsetfillcolor{textcolor}%
\pgftext[x=0.445239in, y=1.300399in, left, base]{\color{textcolor}{\sffamily\fontsize{10.000000}{12.000000}\selectfont\catcode`\^=\active\def^{\ifmmode\sp\else\^{}\fi}\catcode`\%=\active\def
\end{pgfscope}%
\begin{pgfscope}%
\pgfsetbuttcap%
\pgfsetroundjoin%
\definecolor{currentfill}{rgb}{0.000000,0.000000,0.000000}%
\pgfsetfillcolor{currentfill}%
\pgfsetlinewidth{0.803000pt}%
\definecolor{currentstroke}{rgb}{0.000000,0.000000,0.000000}%
\pgfsetstrokecolor{currentstroke}%
\pgfsetdash{}{0pt}%
\pgfsys@defobject{currentmarker}{\pgfqpoint{-0.048611in}{0.000000in}}{\pgfqpoint{-0.000000in}{0.000000in}}{%
\pgfpathmoveto{\pgfqpoint{-0.000000in}{0.000000in}}%
\pgfpathlineto{\pgfqpoint{-0.048611in}{0.000000in}}%
\pgfusepath{stroke,fill}%
}%
\begin{pgfscope}%
\pgfsys@transformshift{0.885827in}{1.921047in}%
\pgfsys@useobject{currentmarker}{}%
\end{pgfscope}%
\end{pgfscope}%
\begin{pgfscope}%
\definecolor{textcolor}{rgb}{0.000000,0.000000,0.000000}%
\pgfsetstrokecolor{textcolor}%
\pgfsetfillcolor{textcolor}%
\pgftext[x=0.500602in, y=1.868286in, left, base]{\color{textcolor}{\sffamily\fontsize{10.000000}{12.000000}\selectfont\catcode`\^=\active\def^{\ifmmode\sp\else\^{}\fi}\catcode`\%=\active\def
\end{pgfscope}%
\begin{pgfscope}%
\pgfsetbuttcap%
\pgfsetroundjoin%
\definecolor{currentfill}{rgb}{0.000000,0.000000,0.000000}%
\pgfsetfillcolor{currentfill}%
\pgfsetlinewidth{0.803000pt}%
\definecolor{currentstroke}{rgb}{0.000000,0.000000,0.000000}%
\pgfsetstrokecolor{currentstroke}%
\pgfsetdash{}{0pt}%
\pgfsys@defobject{currentmarker}{\pgfqpoint{-0.048611in}{0.000000in}}{\pgfqpoint{-0.000000in}{0.000000in}}{%
\pgfpathmoveto{\pgfqpoint{-0.000000in}{0.000000in}}%
\pgfpathlineto{\pgfqpoint{-0.048611in}{0.000000in}}%
\pgfusepath{stroke,fill}%
}%
\begin{pgfscope}%
\pgfsys@transformshift{0.885827in}{2.488934in}%
\pgfsys@useobject{currentmarker}{}%
\end{pgfscope}%
\end{pgfscope}%
\begin{pgfscope}%
\definecolor{textcolor}{rgb}{0.000000,0.000000,0.000000}%
\pgfsetstrokecolor{textcolor}%
\pgfsetfillcolor{textcolor}%
\pgftext[x=0.500602in, y=2.436173in, left, base]{\color{textcolor}{\sffamily\fontsize{10.000000}{12.000000}\selectfont\catcode`\^=\active\def^{\ifmmode\sp\else\^{}\fi}\catcode`\%=\active\def
\end{pgfscope}%
\begin{pgfscope}%
\pgfsetbuttcap%
\pgfsetroundjoin%
\definecolor{currentfill}{rgb}{0.000000,0.000000,0.000000}%
\pgfsetfillcolor{currentfill}%
\pgfsetlinewidth{0.803000pt}%
\definecolor{currentstroke}{rgb}{0.000000,0.000000,0.000000}%
\pgfsetstrokecolor{currentstroke}%
\pgfsetdash{}{0pt}%
\pgfsys@defobject{currentmarker}{\pgfqpoint{-0.048611in}{0.000000in}}{\pgfqpoint{-0.000000in}{0.000000in}}{%
\pgfpathmoveto{\pgfqpoint{-0.000000in}{0.000000in}}%
\pgfpathlineto{\pgfqpoint{-0.048611in}{0.000000in}}%
\pgfusepath{stroke,fill}%
}%
\begin{pgfscope}%
\pgfsys@transformshift{0.885827in}{3.056821in}%
\pgfsys@useobject{currentmarker}{}%
\end{pgfscope}%
\end{pgfscope}%
\begin{pgfscope}%
\definecolor{textcolor}{rgb}{0.000000,0.000000,0.000000}%
\pgfsetstrokecolor{textcolor}%
\pgfsetfillcolor{textcolor}%
\pgftext[x=0.500602in, y=3.004060in, left, base]{\color{textcolor}{\sffamily\fontsize{10.000000}{12.000000}\selectfont\catcode`\^=\active\def^{\ifmmode\sp\else\^{}\fi}\catcode`\%=\active\def
\end{pgfscope}%
\begin{pgfscope}%
\definecolor{textcolor}{rgb}{0.000000,0.000000,0.000000}%
\pgfsetstrokecolor{textcolor}%
\pgfsetfillcolor{textcolor}%
\pgftext[x=0.389684in,y=1.948819in,,bottom,rotate=90.000000]{\color{textcolor}{\sffamily\fontsize{10.000000}{12.000000}\selectfont\catcode`\^=\active\def^{\ifmmode\sp\else\^{}\fi}\catcode`\%=\active\def
\end{pgfscope}%
\begin{pgfscope}%
\pgfpathrectangle{\pgfqpoint{0.885827in}{0.433071in}}{\pgfqpoint{5.492126in}{3.031496in}}%
\pgfusepath{clip}%
\pgfsetrectcap%
\pgfsetroundjoin%
\pgfsetlinewidth{1.505625pt}%
\definecolor{currentstroke}{rgb}{0.000000,0.694118,0.917647}%
\pgfsetstrokecolor{currentstroke}%
\pgfsetdash{}{0pt}%
\pgfpathmoveto{\pgfqpoint{1.266560in}{2.637101in}}%
\pgfpathlineto{\pgfqpoint{2.449225in}{2.811856in}}%
\pgfpathlineto{\pgfqpoint{3.631890in}{2.983778in}}%
\pgfpathlineto{\pgfqpoint{4.814555in}{3.155120in}}%
\pgfpathlineto{\pgfqpoint{5.997220in}{3.326772in}}%
\pgfusepath{stroke}%
\end{pgfscope}%
\begin{pgfscope}%
\pgfpathrectangle{\pgfqpoint{0.885827in}{0.433071in}}{\pgfqpoint{5.492126in}{3.031496in}}%
\pgfusepath{clip}%
\pgfsetbuttcap%
\pgfsetroundjoin%
\definecolor{currentfill}{rgb}{0.000000,0.694118,0.917647}%
\pgfsetfillcolor{currentfill}%
\pgfsetlinewidth{1.003750pt}%
\definecolor{currentstroke}{rgb}{0.000000,0.694118,0.917647}%
\pgfsetstrokecolor{currentstroke}%
\pgfsetdash{}{0pt}%
\pgfsys@defobject{currentmarker}{\pgfqpoint{-0.041667in}{-0.041667in}}{\pgfqpoint{0.041667in}{0.041667in}}{%
\pgfpathmoveto{\pgfqpoint{-0.041667in}{-0.041667in}}%
\pgfpathlineto{\pgfqpoint{0.041667in}{0.041667in}}%
\pgfpathmoveto{\pgfqpoint{-0.041667in}{0.041667in}}%
\pgfpathlineto{\pgfqpoint{0.041667in}{-0.041667in}}%
\pgfusepath{stroke,fill}%
}%
\begin{pgfscope}%
\pgfsys@transformshift{1.266560in}{2.637101in}%
\pgfsys@useobject{currentmarker}{}%
\end{pgfscope}%
\begin{pgfscope}%
\pgfsys@transformshift{2.449225in}{2.811856in}%
\pgfsys@useobject{currentmarker}{}%
\end{pgfscope}%
\begin{pgfscope}%
\pgfsys@transformshift{3.631890in}{2.983778in}%
\pgfsys@useobject{currentmarker}{}%
\end{pgfscope}%
\begin{pgfscope}%
\pgfsys@transformshift{4.814555in}{3.155120in}%
\pgfsys@useobject{currentmarker}{}%
\end{pgfscope}%
\begin{pgfscope}%
\pgfsys@transformshift{5.997220in}{3.326772in}%
\pgfsys@useobject{currentmarker}{}%
\end{pgfscope}%
\end{pgfscope}%
\begin{pgfscope}%
\pgfpathrectangle{\pgfqpoint{0.885827in}{0.433071in}}{\pgfqpoint{5.492126in}{3.031496in}}%
\pgfusepath{clip}%
\pgfsetrectcap%
\pgfsetroundjoin%
\pgfsetlinewidth{1.505625pt}%
\definecolor{currentstroke}{rgb}{0.937255,0.486275,0.000000}%
\pgfsetstrokecolor{currentstroke}%
\pgfsetdash{}{0pt}%
\pgfpathmoveto{\pgfqpoint{1.266560in}{1.628400in}}%
\pgfpathlineto{\pgfqpoint{2.449225in}{1.970499in}}%
\pgfpathlineto{\pgfqpoint{3.631890in}{2.312448in}}%
\pgfpathlineto{\pgfqpoint{4.814555in}{2.654493in}}%
\pgfpathlineto{\pgfqpoint{5.997220in}{2.996963in}}%
\pgfusepath{stroke}%
\end{pgfscope}%
\begin{pgfscope}%
\pgfpathrectangle{\pgfqpoint{0.885827in}{0.433071in}}{\pgfqpoint{5.492126in}{3.031496in}}%
\pgfusepath{clip}%
\pgfsetbuttcap%
\pgfsetroundjoin%
\definecolor{currentfill}{rgb}{0.937255,0.486275,0.000000}%
\pgfsetfillcolor{currentfill}%
\pgfsetlinewidth{1.003750pt}%
\definecolor{currentstroke}{rgb}{0.937255,0.486275,0.000000}%
\pgfsetstrokecolor{currentstroke}%
\pgfsetdash{}{0pt}%
\pgfsys@defobject{currentmarker}{\pgfqpoint{-0.041667in}{-0.041667in}}{\pgfqpoint{0.041667in}{0.041667in}}{%
\pgfpathmoveto{\pgfqpoint{-0.041667in}{-0.041667in}}%
\pgfpathlineto{\pgfqpoint{0.041667in}{0.041667in}}%
\pgfpathmoveto{\pgfqpoint{-0.041667in}{0.041667in}}%
\pgfpathlineto{\pgfqpoint{0.041667in}{-0.041667in}}%
\pgfusepath{stroke,fill}%
}%
\begin{pgfscope}%
\pgfsys@transformshift{1.266560in}{1.628400in}%
\pgfsys@useobject{currentmarker}{}%
\end{pgfscope}%
\begin{pgfscope}%
\pgfsys@transformshift{2.449225in}{1.970499in}%
\pgfsys@useobject{currentmarker}{}%
\end{pgfscope}%
\begin{pgfscope}%
\pgfsys@transformshift{3.631890in}{2.312448in}%
\pgfsys@useobject{currentmarker}{}%
\end{pgfscope}%
\begin{pgfscope}%
\pgfsys@transformshift{4.814555in}{2.654493in}%
\pgfsys@useobject{currentmarker}{}%
\end{pgfscope}%
\begin{pgfscope}%
\pgfsys@transformshift{5.997220in}{2.996963in}%
\pgfsys@useobject{currentmarker}{}%
\end{pgfscope}%
\end{pgfscope}%
\begin{pgfscope}%
\pgfpathrectangle{\pgfqpoint{0.885827in}{0.433071in}}{\pgfqpoint{5.492126in}{3.031496in}}%
\pgfusepath{clip}%
\pgfsetrectcap%
\pgfsetroundjoin%
\pgfsetlinewidth{1.505625pt}%
\definecolor{currentstroke}{rgb}{0.000000,0.423529,0.400000}%
\pgfsetstrokecolor{currentstroke}%
\pgfsetdash{}{0pt}%
\pgfpathmoveto{\pgfqpoint{1.266560in}{1.605204in}}%
\pgfpathlineto{\pgfqpoint{2.449225in}{1.779949in}}%
\pgfpathlineto{\pgfqpoint{3.631890in}{1.951828in}}%
\pgfpathlineto{\pgfqpoint{4.814555in}{2.122999in}}%
\pgfpathlineto{\pgfqpoint{5.997220in}{2.293962in}}%
\pgfusepath{stroke}%
\end{pgfscope}%
\begin{pgfscope}%
\pgfpathrectangle{\pgfqpoint{0.885827in}{0.433071in}}{\pgfqpoint{5.492126in}{3.031496in}}%
\pgfusepath{clip}%
\pgfsetbuttcap%
\pgfsetroundjoin%
\definecolor{currentfill}{rgb}{0.000000,0.423529,0.400000}%
\pgfsetfillcolor{currentfill}%
\pgfsetlinewidth{1.003750pt}%
\definecolor{currentstroke}{rgb}{0.000000,0.423529,0.400000}%
\pgfsetstrokecolor{currentstroke}%
\pgfsetdash{}{0pt}%
\pgfsys@defobject{currentmarker}{\pgfqpoint{-0.041667in}{-0.041667in}}{\pgfqpoint{0.041667in}{0.041667in}}{%
\pgfpathmoveto{\pgfqpoint{-0.041667in}{-0.041667in}}%
\pgfpathlineto{\pgfqpoint{0.041667in}{0.041667in}}%
\pgfpathmoveto{\pgfqpoint{-0.041667in}{0.041667in}}%
\pgfpathlineto{\pgfqpoint{0.041667in}{-0.041667in}}%
\pgfusepath{stroke,fill}%
}%
\begin{pgfscope}%
\pgfsys@transformshift{1.266560in}{1.605204in}%
\pgfsys@useobject{currentmarker}{}%
\end{pgfscope}%
\begin{pgfscope}%
\pgfsys@transformshift{2.449225in}{1.779949in}%
\pgfsys@useobject{currentmarker}{}%
\end{pgfscope}%
\begin{pgfscope}%
\pgfsys@transformshift{3.631890in}{1.951828in}%
\pgfsys@useobject{currentmarker}{}%
\end{pgfscope}%
\begin{pgfscope}%
\pgfsys@transformshift{4.814555in}{2.122999in}%
\pgfsys@useobject{currentmarker}{}%
\end{pgfscope}%
\begin{pgfscope}%
\pgfsys@transformshift{5.997220in}{2.293962in}%
\pgfsys@useobject{currentmarker}{}%
\end{pgfscope}%
\end{pgfscope}%
\begin{pgfscope}%
\pgfpathrectangle{\pgfqpoint{0.885827in}{0.433071in}}{\pgfqpoint{5.492126in}{3.031496in}}%
\pgfusepath{clip}%
\pgfsetrectcap%
\pgfsetroundjoin%
\pgfsetlinewidth{1.505625pt}%
\definecolor{currentstroke}{rgb}{0.776471,0.827451,0.145098}%
\pgfsetstrokecolor{currentstroke}%
\pgfsetdash{}{0pt}%
\pgfpathmoveto{\pgfqpoint{1.266560in}{0.570866in}}%
\pgfpathlineto{\pgfqpoint{2.449225in}{0.609126in}}%
\pgfpathlineto{\pgfqpoint{3.631890in}{0.879107in}}%
\pgfpathlineto{\pgfqpoint{4.814555in}{1.220266in}}%
\pgfpathlineto{\pgfqpoint{5.997220in}{1.562053in}}%
\pgfusepath{stroke}%
\end{pgfscope}%
\begin{pgfscope}%
\pgfpathrectangle{\pgfqpoint{0.885827in}{0.433071in}}{\pgfqpoint{5.492126in}{3.031496in}}%
\pgfusepath{clip}%
\pgfsetbuttcap%
\pgfsetroundjoin%
\definecolor{currentfill}{rgb}{0.776471,0.827451,0.145098}%
\pgfsetfillcolor{currentfill}%
\pgfsetlinewidth{1.003750pt}%
\definecolor{currentstroke}{rgb}{0.776471,0.827451,0.145098}%
\pgfsetstrokecolor{currentstroke}%
\pgfsetdash{}{0pt}%
\pgfsys@defobject{currentmarker}{\pgfqpoint{-0.041667in}{-0.041667in}}{\pgfqpoint{0.041667in}{0.041667in}}{%
\pgfpathmoveto{\pgfqpoint{-0.041667in}{-0.041667in}}%
\pgfpathlineto{\pgfqpoint{0.041667in}{0.041667in}}%
\pgfpathmoveto{\pgfqpoint{-0.041667in}{0.041667in}}%
\pgfpathlineto{\pgfqpoint{0.041667in}{-0.041667in}}%
\pgfusepath{stroke,fill}%
}%
\begin{pgfscope}%
\pgfsys@transformshift{1.266560in}{0.570866in}%
\pgfsys@useobject{currentmarker}{}%
\end{pgfscope}%
\begin{pgfscope}%
\pgfsys@transformshift{2.449225in}{0.609126in}%
\pgfsys@useobject{currentmarker}{}%
\end{pgfscope}%
\begin{pgfscope}%
\pgfsys@transformshift{3.631890in}{0.879107in}%
\pgfsys@useobject{currentmarker}{}%
\end{pgfscope}%
\begin{pgfscope}%
\pgfsys@transformshift{4.814555in}{1.220266in}%
\pgfsys@useobject{currentmarker}{}%
\end{pgfscope}%
\begin{pgfscope}%
\pgfsys@transformshift{5.997220in}{1.562053in}%
\pgfsys@useobject{currentmarker}{}%
\end{pgfscope}%
\end{pgfscope}%
\begin{pgfscope}%
\pgfsetrectcap%
\pgfsetmiterjoin%
\pgfsetlinewidth{0.803000pt}%
\definecolor{currentstroke}{rgb}{0.000000,0.000000,0.000000}%
\pgfsetstrokecolor{currentstroke}%
\pgfsetdash{}{0pt}%
\pgfpathmoveto{\pgfqpoint{0.885827in}{0.433071in}}%
\pgfpathlineto{\pgfqpoint{0.885827in}{3.464567in}}%
\pgfusepath{stroke}%
\end{pgfscope}%
\begin{pgfscope}%
\pgfsetrectcap%
\pgfsetmiterjoin%
\pgfsetlinewidth{0.803000pt}%
\definecolor{currentstroke}{rgb}{0.000000,0.000000,0.000000}%
\pgfsetstrokecolor{currentstroke}%
\pgfsetdash{}{0pt}%
\pgfpathmoveto{\pgfqpoint{6.377953in}{0.433071in}}%
\pgfpathlineto{\pgfqpoint{6.377953in}{3.464567in}}%
\pgfusepath{stroke}%
\end{pgfscope}%
\begin{pgfscope}%
\pgfsetrectcap%
\pgfsetmiterjoin%
\pgfsetlinewidth{0.803000pt}%
\definecolor{currentstroke}{rgb}{0.000000,0.000000,0.000000}%
\pgfsetstrokecolor{currentstroke}%
\pgfsetdash{}{0pt}%
\pgfpathmoveto{\pgfqpoint{0.885827in}{0.433071in}}%
\pgfpathlineto{\pgfqpoint{6.377953in}{0.433071in}}%
\pgfusepath{stroke}%
\end{pgfscope}%
\begin{pgfscope}%
\pgfsetrectcap%
\pgfsetmiterjoin%
\pgfsetlinewidth{0.803000pt}%
\definecolor{currentstroke}{rgb}{0.000000,0.000000,0.000000}%
\pgfsetstrokecolor{currentstroke}%
\pgfsetdash{}{0pt}%
\pgfpathmoveto{\pgfqpoint{0.885827in}{3.464567in}}%
\pgfpathlineto{\pgfqpoint{6.377953in}{3.464567in}}%
\pgfusepath{stroke}%
\end{pgfscope}%
\begin{pgfscope}%
\pgfsetbuttcap%
\pgfsetmiterjoin%
\definecolor{currentfill}{rgb}{1.000000,1.000000,1.000000}%
\pgfsetfillcolor{currentfill}%
\pgfsetfillopacity{0.800000}%
\pgfsetlinewidth{1.003750pt}%
\definecolor{currentstroke}{rgb}{0.800000,0.800000,0.800000}%
\pgfsetstrokecolor{currentstroke}%
\pgfsetstrokeopacity{0.800000}%
\pgfsetdash{}{0pt}%
\pgfpathmoveto{\pgfqpoint{0.983049in}{2.538027in}}%
\pgfpathlineto{\pgfqpoint{3.019697in}{2.538027in}}%
\pgfpathquadraticcurveto{\pgfqpoint{3.047475in}{2.538027in}}{\pgfqpoint{3.047475in}{2.565805in}}%
\pgfpathlineto{\pgfqpoint{3.047475in}{3.367345in}}%
\pgfpathquadraticcurveto{\pgfqpoint{3.047475in}{3.395122in}}{\pgfqpoint{3.019697in}{3.395122in}}%
\pgfpathlineto{\pgfqpoint{0.983049in}{3.395122in}}%
\pgfpathquadraticcurveto{\pgfqpoint{0.955271in}{3.395122in}}{\pgfqpoint{0.955271in}{3.367345in}}%
\pgfpathlineto{\pgfqpoint{0.955271in}{2.565805in}}%
\pgfpathquadraticcurveto{\pgfqpoint{0.955271in}{2.538027in}}{\pgfqpoint{0.983049in}{2.538027in}}%
\pgfpathlineto{\pgfqpoint{0.983049in}{2.538027in}}%
\pgfpathclose%
\pgfusepath{stroke,fill}%
\end{pgfscope}%
\begin{pgfscope}%
\pgfsetrectcap%
\pgfsetroundjoin%
\pgfsetlinewidth{1.505625pt}%
\definecolor{currentstroke}{rgb}{0.000000,0.694118,0.917647}%
\pgfsetstrokecolor{currentstroke}%
\pgfsetdash{}{0pt}%
\pgfpathmoveto{\pgfqpoint{1.010827in}{3.282655in}}%
\pgfpathlineto{\pgfqpoint{1.149716in}{3.282655in}}%
\pgfpathlineto{\pgfqpoint{1.288605in}{3.282655in}}%
\pgfusepath{stroke}%
\end{pgfscope}%
\begin{pgfscope}%
\pgfsetbuttcap%
\pgfsetroundjoin%
\definecolor{currentfill}{rgb}{0.000000,0.694118,0.917647}%
\pgfsetfillcolor{currentfill}%
\pgfsetlinewidth{1.003750pt}%
\definecolor{currentstroke}{rgb}{0.000000,0.694118,0.917647}%
\pgfsetstrokecolor{currentstroke}%
\pgfsetdash{}{0pt}%
\pgfsys@defobject{currentmarker}{\pgfqpoint{-0.041667in}{-0.041667in}}{\pgfqpoint{0.041667in}{0.041667in}}{%
\pgfpathmoveto{\pgfqpoint{-0.041667in}{-0.041667in}}%
\pgfpathlineto{\pgfqpoint{0.041667in}{0.041667in}}%
\pgfpathmoveto{\pgfqpoint{-0.041667in}{0.041667in}}%
\pgfpathlineto{\pgfqpoint{0.041667in}{-0.041667in}}%
\pgfusepath{stroke,fill}%
}%
\begin{pgfscope}%
\pgfsys@transformshift{1.149716in}{3.282655in}%
\pgfsys@useobject{currentmarker}{}%
\end{pgfscope}%
\end{pgfscope}%
\begin{pgfscope}%
\definecolor{textcolor}{rgb}{0.000000,0.000000,0.000000}%
\pgfsetstrokecolor{textcolor}%
\pgfsetfillcolor{textcolor}%
\pgftext[x=1.399716in,y=3.234044in,left,base]{\color{textcolor}{\sffamily\fontsize{10.000000}{12.000000}\selectfont\catcode`\^=\active\def^{\ifmmode\sp\else\^{}\fi}\catcode`\%=\active\def
\end{pgfscope}%
\begin{pgfscope}%
\pgfsetrectcap%
\pgfsetroundjoin%
\pgfsetlinewidth{1.505625pt}%
\definecolor{currentstroke}{rgb}{0.937255,0.486275,0.000000}%
\pgfsetstrokecolor{currentstroke}%
\pgfsetdash{}{0pt}%
\pgfpathmoveto{\pgfqpoint{1.010827in}{3.078798in}}%
\pgfpathlineto{\pgfqpoint{1.149716in}{3.078798in}}%
\pgfpathlineto{\pgfqpoint{1.288605in}{3.078798in}}%
\pgfusepath{stroke}%
\end{pgfscope}%
\begin{pgfscope}%
\pgfsetbuttcap%
\pgfsetroundjoin%
\definecolor{currentfill}{rgb}{0.937255,0.486275,0.000000}%
\pgfsetfillcolor{currentfill}%
\pgfsetlinewidth{1.003750pt}%
\definecolor{currentstroke}{rgb}{0.937255,0.486275,0.000000}%
\pgfsetstrokecolor{currentstroke}%
\pgfsetdash{}{0pt}%
\pgfsys@defobject{currentmarker}{\pgfqpoint{-0.041667in}{-0.041667in}}{\pgfqpoint{0.041667in}{0.041667in}}{%
\pgfpathmoveto{\pgfqpoint{-0.041667in}{-0.041667in}}%
\pgfpathlineto{\pgfqpoint{0.041667in}{0.041667in}}%
\pgfpathmoveto{\pgfqpoint{-0.041667in}{0.041667in}}%
\pgfpathlineto{\pgfqpoint{0.041667in}{-0.041667in}}%
\pgfusepath{stroke,fill}%
}%
\begin{pgfscope}%
\pgfsys@transformshift{1.149716in}{3.078798in}%
\pgfsys@useobject{currentmarker}{}%
\end{pgfscope}%
\end{pgfscope}%
\begin{pgfscope}%
\definecolor{textcolor}{rgb}{0.000000,0.000000,0.000000}%
\pgfsetstrokecolor{textcolor}%
\pgfsetfillcolor{textcolor}%
\pgftext[x=1.399716in,y=3.030187in,left,base]{\color{textcolor}{\sffamily\fontsize{10.000000}{12.000000}\selectfont\catcode`\^=\active\def^{\ifmmode\sp\else\^{}\fi}\catcode`\%=\active\def
\end{pgfscope}%
\begin{pgfscope}%
\pgfsetrectcap%
\pgfsetroundjoin%
\pgfsetlinewidth{1.505625pt}%
\definecolor{currentstroke}{rgb}{0.000000,0.423529,0.400000}%
\pgfsetstrokecolor{currentstroke}%
\pgfsetdash{}{0pt}%
\pgfpathmoveto{\pgfqpoint{1.010827in}{2.874941in}}%
\pgfpathlineto{\pgfqpoint{1.149716in}{2.874941in}}%
\pgfpathlineto{\pgfqpoint{1.288605in}{2.874941in}}%
\pgfusepath{stroke}%
\end{pgfscope}%
\begin{pgfscope}%
\pgfsetbuttcap%
\pgfsetroundjoin%
\definecolor{currentfill}{rgb}{0.000000,0.423529,0.400000}%
\pgfsetfillcolor{currentfill}%
\pgfsetlinewidth{1.003750pt}%
\definecolor{currentstroke}{rgb}{0.000000,0.423529,0.400000}%
\pgfsetstrokecolor{currentstroke}%
\pgfsetdash{}{0pt}%
\pgfsys@defobject{currentmarker}{\pgfqpoint{-0.041667in}{-0.041667in}}{\pgfqpoint{0.041667in}{0.041667in}}{%
\pgfpathmoveto{\pgfqpoint{-0.041667in}{-0.041667in}}%
\pgfpathlineto{\pgfqpoint{0.041667in}{0.041667in}}%
\pgfpathmoveto{\pgfqpoint{-0.041667in}{0.041667in}}%
\pgfpathlineto{\pgfqpoint{0.041667in}{-0.041667in}}%
\pgfusepath{stroke,fill}%
}%
\begin{pgfscope}%
\pgfsys@transformshift{1.149716in}{2.874941in}%
\pgfsys@useobject{currentmarker}{}%
\end{pgfscope}%
\end{pgfscope}%
\begin{pgfscope}%
\definecolor{textcolor}{rgb}{0.000000,0.000000,0.000000}%
\pgfsetstrokecolor{textcolor}%
\pgfsetfillcolor{textcolor}%
\pgftext[x=1.399716in,y=2.826329in,left,base]{\color{textcolor}{\sffamily\fontsize{10.000000}{12.000000}\selectfont\catcode`\^=\active\def^{\ifmmode\sp\else\^{}\fi}\catcode`\%=\active\def
\end{pgfscope}%
\begin{pgfscope}%
\pgfsetrectcap%
\pgfsetroundjoin%
\pgfsetlinewidth{1.505625pt}%
\definecolor{currentstroke}{rgb}{0.776471,0.827451,0.145098}%
\pgfsetstrokecolor{currentstroke}%
\pgfsetdash{}{0pt}%
\pgfpathmoveto{\pgfqpoint{1.010827in}{2.671083in}}%
\pgfpathlineto{\pgfqpoint{1.149716in}{2.671083in}}%
\pgfpathlineto{\pgfqpoint{1.288605in}{2.671083in}}%
\pgfusepath{stroke}%
\end{pgfscope}%
\begin{pgfscope}%
\pgfsetbuttcap%
\pgfsetroundjoin%
\definecolor{currentfill}{rgb}{0.776471,0.827451,0.145098}%
\pgfsetfillcolor{currentfill}%
\pgfsetlinewidth{1.003750pt}%
\definecolor{currentstroke}{rgb}{0.776471,0.827451,0.145098}%
\pgfsetstrokecolor{currentstroke}%
\pgfsetdash{}{0pt}%
\pgfsys@defobject{currentmarker}{\pgfqpoint{-0.041667in}{-0.041667in}}{\pgfqpoint{0.041667in}{0.041667in}}{%
\pgfpathmoveto{\pgfqpoint{-0.041667in}{-0.041667in}}%
\pgfpathlineto{\pgfqpoint{0.041667in}{0.041667in}}%
\pgfpathmoveto{\pgfqpoint{-0.041667in}{0.041667in}}%
\pgfpathlineto{\pgfqpoint{0.041667in}{-0.041667in}}%
\pgfusepath{stroke,fill}%
}%
\begin{pgfscope}%
\pgfsys@transformshift{1.149716in}{2.671083in}%
\pgfsys@useobject{currentmarker}{}%
\end{pgfscope}%
\end{pgfscope}%
\begin{pgfscope}%
\definecolor{textcolor}{rgb}{0.000000,0.000000,0.000000}%
\pgfsetstrokecolor{textcolor}%
\pgfsetfillcolor{textcolor}%
\pgftext[x=1.399716in,y=2.622472in,left,base]{\color{textcolor}{\sffamily\fontsize{10.000000}{12.000000}\selectfont\catcode`\^=\active\def^{\ifmmode\sp\else\^{}\fi}\catcode`\%=\active\def
\end{pgfscope}%
\end{pgfpicture}%
\makeatother%
\endgroup%

%% file: Paper.bbl
\begin{thebibliography}{10}

\bibitem{Bernier2020}
{\sc J.~Bernier, F.~Casas, and N.~Crouseilles}, {\em Splitting {Methods} for {Rotations}: {Application} to {Vlasov} {Equations}}, SIAM J. Sci. Comput., 42 (2020), pp.~A666--A697, \url{https://doi.org/10.1137/19M1273918}.

\bibitem{Besse2008}
{\sc N.~Besse and M.~Mehrenberger}, {\em Convergence of classes of high-order semi-{Lagrangian} schemes for the {Vlasov}--{Poisson} system}, Math. Comp., 77 (2008), pp.~93--123, \url{https://doi.org/10.1090/S0025-5718-07-01912-6}.

\bibitem{brambilla_kinetic_1998}
{\sc M.~Brambilla}, {\em Kinetic theory of plasma waves: homogeneous plasmas}, The International series of monographs on physics, Clarendon Press, Oxford ; New York, 1998.

\bibitem{BSL6D}
{\sc {BSL6D Development Team}}, {\em {BSL6D}}, 2015, \url{https://gitlab.mpcdf.mpg.de/bsl6d/bsl6d} (accessed 2024/05/02).

\bibitem{Chen2016}
{\sc F.~F. Chen}, {\em Introduction to {Plasma} {Physics} and {Controlled} {Fusion}}, Springer, Cham, 3rd ed. 2016~ed., 2016.

\bibitem{Cheng1976}
{\sc C.~Cheng and G.~Knorr}, {\em The integration of the vlasov equation in configuration space}, Journal of Computational Physics, 22 (1976), pp.~330--351, \url{https://doi.org/10.1016/0021-9991(76)90053-X}.

\bibitem{donea2004}
{\sc J.~Donea, A.~Huerta, J.-P. Ponthot, and A.~Rodríguez-Ferran}, {\em Arbitrary Lagrangian-Eulerian Methods}, John Wiley {\&} Sons, Ltd, Nov. 2004, ch.~14, pp.~413--437, \url{https://doi.org/10.1002/0470091355.ecm009}.

\bibitem{Einkemmer2014}
{\sc L.~Einkemmer and A.~Ostermann}, {\em Convergence {Analysis} of {Strang} {Splitting} for {Vlasov}-{Type} {Equations}}, SIAM J. Numer. Anal., 52 (2014), pp.~140--155, \url{https://doi.org/10.1137/130918599}.

\bibitem{Falcone2014}
{\sc M.~Falcone and R.~Ferretti}, {\em Semi-{Lagrangian} approximation schemes for linear and {Hamilton}-{Jacobi} equations}, {OT} / {SIAM}, {Society} of {Industrial} and {Applied} {Mathematics}, SIAM, Philadelphia, Pa, 2014.

\bibitem{Hairer2006}
{\sc E.~Hairer, C.~Lubich, and G.~Wanner}, {\em Geometric numerical integration: structure-preserving algorithms for ordinary differential equations}, Springer series in computational mathematics, Springer, Berlin ; New York, 2nd ed~ed., 2006.

\bibitem{Huang2011}
{\sc W.~Huang and R.~D. Russell}, {\em Adaptive moving mesh methods}, Applied mathematical sciences ({Springer}-{Verlag} {New} {York} {Inc}.), Springer, New York, NY, 2011.

\bibitem{Kormann2019}
{\sc K.~Kormann, K.~Reuter, and M.~Rampp}, {\em A massively parallel semi-{Lagrangian} solver for the six-dimensional {Vlasov}–{Poisson} equation}, The International Journal of High Performance Computing Applications, 33 (2019), pp.~924--947, \url{https://doi.org/10.1177/1094342019834644}.

\bibitem{Kraus2017}
{\sc M.~Kraus, K.~Kormann, P.~J. Morrison, and E.~Sonnendrücker}, {\em {GEMPIC}: geometric electromagnetic particle-in-cell methods}, J. Plasma Phys., 83 (2017), pp.~905830401--905830451, \url{https://doi.org/10.1017/S002237781700040X}.

\bibitem{Mukhamet2023}
{\sc T.~Mukhamet}, {\em An arbitrary {L}agrangian {E}ulerian discontinuous galerkin method for {V}lasov equation with a strong magnetic field}, 2023.

\bibitem{RaethThesis2023}
{\sc M.~Raeth}, {\em Beyond gyrokinetic theory: Excitation of high-frequency turbulence in 6d {V}lasov simulations of magnetized plasmas with steep temperature and density gradients}, 2023, \url{https://mediatum.ub.tum.de/node?id=1703830&change_language=en}.

\bibitem{RaethPhysRevLett24}
{\sc M.~Raeth and K.~Hallatschek}, {\em High frequency non-gyrokinetic turbulence at tokamak edge parameters}, Phys. Rev. Lett.,  (2023), \url{https://doi.org/10.48550/arXiv.2310.15981}.

\bibitem{RaethPhysPlasmas24}
{\sc M.~Raeth, K.~Hallatschek, and K.~Kormann}, {\em {Simulation of ion temperature gradient driven modes with 6D kinetic Vlasov code}}, Physics of Plasmas, 31 (2024), p.~042101, \url{https://doi.org/10.1063/5.0197970}.

\bibitem{Schild2024}
{\sc N.~Schild, M.~Räth, S.~Eibl, K.~Hallatschek, and K.~Kormann}, {\em A performance portable implementation of the semi-{Lagrangian} algorithm in six dimensions}, Computer Physics Communications, 295 (2024), p.~108973, \url{https://doi.org/10.1016/j.cpc.2023.108973}.

\bibitem{Schmitz2006}
{\sc H.~Schmitz and R.~Grauer}, {\em Comparison of time splitting and backsubstitution methods for integrating {V}lasov's equation with magnetic fields}, Computer Physics Communications, 175 (2006), pp.~86--92, \url{https://doi.org/10.1016/j.cpc.2006.02.007}.

\end{thebibliography}
